\documentclass[12pt,fleqn,reqno]{amsbook}

\usepackage{culmus} 
\makeatletter
\def\@rllanguagename{hebrew}
\def\eqalign#1{%
 \null\,\vcenter{\openup\jot\m@th
  \ialign{\strut\hfil$\displaystyle{##}$&$\displaystyle{{}##}$\hfil
      \crcr#1\crcr}}\,}
\makeatother

\usepackage{amstext}
\usepackage{amssymb}
\usepackage{amsmath}
\usepackage{amsthm}
\usepackage{amsopn}

\usepackage{graphicx}

\theoremstyle{plain}
\newtheorem{thm}{Theorem}

\newtheorem{lem}[thm]{Lemma}
\newtheorem{problem}[thm]{Problem}

\theoremstyle{definition}

\newtheorem{note}[thm]{Remark}

\DeclareMathOperator{\CMD}{gcd}
\DeclareMathOperator{\CMM}{lcm}

\DeclareMathOperator{\im}{Im}
\DeclareMathOperator{\re}{Re}

\DeclareMathOperator{\Log}{Log}
\DeclareMathOperator{\Arg}{Arg}
\DeclareMathOperator{\Ind}{Ind}
\DeclareMathOperator{\PR}{P\mathbb R}
\DeclareMathOperator{\sgn}{sign}
\DeclareMathOperator{\ran}{rank}
\DeclareMathOperator{\Pos}{pos}
\DeclareMathOperator{\Neg}{neg}
\DeclareMathOperator{\codim}{codim}
\DeclareMathOperator{\Res}{Res}

\let\originalR\R
\newcommand{\newR}{\mathbb R}
\let\R\newR

\newcommand{\IndI}{\Ind_\infty}
\newcommand{\IndC}{\Ind_{-\infty}^{+\infty}}
\newcommand{\IndPR}{\Ind_{\PR}}

\newcommand\SL{\pmb{S}\pmb{L}(2,\mathbb{R})}
\newcommand\Routh{{\tt Routh\,}}

\let\phi\varphi
\let\epsilon\varepsilon
\let\kappa\varkappa
\let\leq\leqslant

\newcommand\tvr{}

\newcommand\resume{\small\sc}
\newcommand{\spec}[2]{\hspace{#1 mm}\tt #2}

\begin{document}

\makeatletter
\renewcommand{\@oddhead}{\hfil\thepage\hfil }
\renewcommand{\@evenhead}{\hfil\thepage\hfil }
\makeatother

\thispagestyle{empty}


\title{Lectures on the Routh-Hurwitz problem}
\author{Yury S. Barkovsky  \\  \small Department of Mathematics, Mechanics and Computer Science \\
 \small South Federal University,
  Rostov-on-Don, Russia \\[0.6cm]
 \small translated from the Russian by
 Olga V. Holtz and Mikhail Yu. Tyaglov  }

\date{\today}
\maketitle

\tableofcontents

\chapter{Basic theory}
\vskip 3cm
\begin{figure}[ht]
\includegraphics[width=160mm, keepaspectratio=yes]{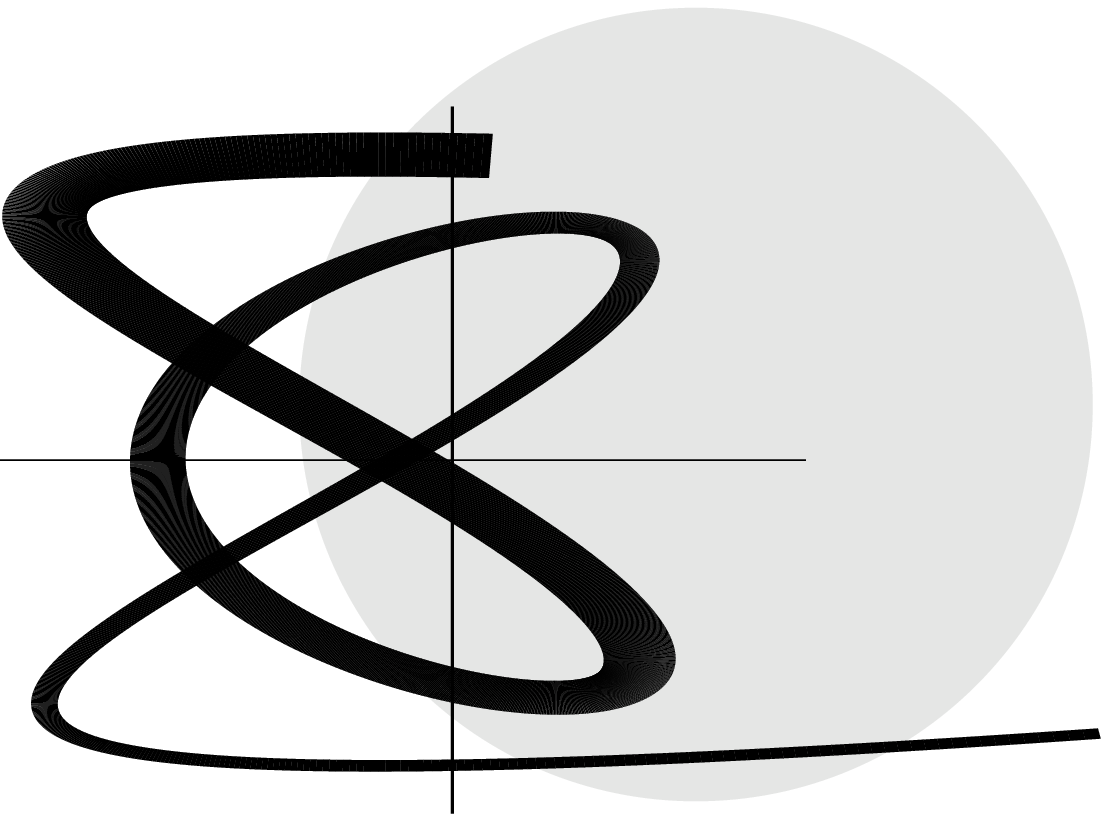}
\end{figure}
\vskip 2cm
\begin{quote}
''Well'', said Owl, ''the customary procedure in
such cases is as follows.''
\end{quote}

\begin{flushright}
A. A. Milne, {\sc Winnie-the-Pooh.}
\end{flushright}
\bigskip

\newpage
\section*{Introduction}

In the vicinity of the zero fixed point%
\footnote{Lit.: point of equilibrium [translators' remark].}
 of the differential equation
\begin{equation}\label{eq1}
    \dot{x}=Ax+o(x)\qquad(x\in{\R}^n),
\end{equation}
the behaviour of other solutions is mostly determined%
\footnote{In the so-called {\it critical case,\/} the nonlinear term
$o(x)$ <<gets voting rights>> and influences the behaviour of
solutions in an arbitrarily small vicinity of the fixed point.} by
the location of the eigenvalues of the matrix~$A$ with respect to the
imaginary axis. For instance, if the spectrum of the matrix $A$ lies
in the  open left half-plane of the complex plane, then the fixed point
is asymptotically stable. It is unstable if the matrix~$A$
has {\it at least\/} one eigenvalue with positive real part. If
we know that the matrix $A$ has no purely imaginary eigenvalues,
then the local structure of the phase portrait of the system~\eqref{eq1}
can be determined by the number of eigenvalues of the matrix $A$ in the
left and right half-planes of the complex plane. For all these reasons,
matrix theory methods and techniques for answering such questions
are of great interest in stability theory. Since the spectrum
of the matrix $A$ is the set of all roots of its characteristic
polynomial~$p(\lambda)\equiv\det(\lambda I-A)$, then the same
questions are to be answered in the theory of polynomials, too.%
\footnote {The transition from a matrix to its characteristic polynomial
is far from being harmless. Firstly, matrix properties that may have
an influence on its spectrum may be lost or hidden as a result.
For example, it is easy to establish that all eigenvalues of a
symmetric matrix must be real, but it is more difficult to understand
how this property affects  the coefficients of its characteristic
polynomial. Secondly, the characteristic polynomial is useful only
for general theoretical questions and does not easily submit to
numerical computations and analytical derivations.}

A polynomial $p$ is called {\it stable\/} if all its roots lie in the open
left half-plane. The {\it Routh-Hurwitz problem\/} consists in finding
conditions of polynomial stability%
\footnote
{Both terms became customary but they are rather unsatisfactory. The former
is not good because there is a <<stable>> fixed point of a system of
differential equations not of a polynomial. The latter is also
bad, since the first person who posed the  <<Routh-Hurwitz problem>>
and who obtained fundamental results in this area was in fact C.\,Hermite.
Here is the chronology of works:
C.\,Hermite --- 1856, E.\,J.\,Routh --- 1877, A.\,Hurwitz --- 1895.}
and, generally, in the study of those properties of the polynomial that
are  in some way connected with location of its roots with respect to
the imaginary axis.

The Hurwitz criterion is traditionally viewed as the main result on stable
polynomials. It will be discussed in \S\,\ref{sec6}. The practical use of this
theorem is usually limited, in the context of direct computations, to
polynomials of low degrees ($3$rd, $4$th, or $5$th).
In fact, the Hurwitz criterion is only one of the facts of a
compact algebraic theory -- a theory that contains other practically useful
results, that is related to important and interesting chapters of
algebra and analysis, and, finally, that is beautiful. The last point
is important. It is most difficult to master the art of producing
mathematical results, of posing and solving problems, but
it is an art worth learning%
\footnote{The reader is referred to textbooks \cite{8}--\cite{10}
and problem book \cite{11a,11b}.}.  Unexpected ideas, subtle arguments
are seldom fruits of pure imagination; more often, they are results
of observation, perseverance, and good taste of their author.
That kind of experience comes with learning things that are worth emulating.
The theory of stable polynomials provides a great sample of this
kind. Within this theory, everyday mathematical notions and ideas
interact, reshape themselves, and bring about new realms of
possible applications. Taking these didactic ideas to heart, the
author did not intend  to simply give a standard list of facts,
but instead to show the development of this mathematical theory,
so that  the reader may became a participant in its re-creation.
Acknowledging that the interest of some student readers may be quite
pragmatic, the author at the same time tried to separate the basic
material, which one ought to learn in any case and which is
presented very tersely,  from discussions and additional
points made in remarks, problems, footnotes etc. Incidentally,
one can learn to apply the Hurwitz theorem by solving the following
fairly typical problem. \vskip 0.4cm

\noindent{\bf Problem.} {\it
Find all fixed points of the {\it Lorentz system}%
\footnote{This system is related to one of the classical
hydrodynamical problems, viz., the onset problem for convectional
motion in a fluid horizontal layer heated from below. The Lorentz system
is interesting due to the fact that its trajectories have very complex
behaviour for certain values of parameters.}
\begin{equation}\label{eq2}
    \dot{X}=\sigma Y-\sigma X,\quad
    \dot{Y}=rX-Y-XZ,\quad
    \dot{Z}=XY-bZ.
\end{equation}
($\sigma$, $r$, $b$ are positive parameters).
Get to know the statement of Hurwitz' theorem in \S\,\ref{sec6} and
apply it to investigate the stability of the fixed points found.}

A solution to this problem is given in the Appendix. It is however recommended
that the reader obtains this solution on her/his own or at least tries to
do so.

\section*{Acknowledgment}   \noindent
The author thanks Mikhail Tyaglov, who pointed out a number of typos 
in the first version of these notes; those typos are corrected in the
present version. The author is also grateful to Olga Holtz for her 
active interest in these notes and assistance with their publication.

\section{Stodola condition}\label{sec1}

One of the most basic but rather useful facts on stable polynomials is 
contained in the following theorem, which is usually attributed to the
Slovak engineer A. Stodola (1893).
\begin{thm}[Stodola]\label{thm1.1}
If a polynomial with real coefficients is stable, then all its
coefficients are of the same sign.
\end{thm}
\begin{proof}
The roots of a real polynomial are symmetric with respect to the real axis.
Let
\begin{equation*}
p(z)=a_0\prod\limits_j\left(z-\lambda_j\right)\cdot
        \prod\limits_k\left(z-\alpha_k-i\beta_k\right)
                      \left(z-\alpha_k+i\beta_k\right),
\end{equation*}
where $\lambda_j$ are the real and $\alpha_k\pm i\beta_k$ are the nonreal 
roots of the polynomial $p$ (note that $\lambda_j,\alpha_k<0$). 
Since the binomials~$z-\lambda_j$ and the trinomials~$z^2-2\alpha_kz+(\alpha_k^2+\beta_k^2)$ have positive coefficients, their product has the same property.
\end{proof}

The role of Theorem~\ref{thm1.1} is quite clear: it provides a very easily 
verified {\it necessary\/} condition of polynomial stability. It  cannot be 
reversed, except for very low degrees:
\begin{problem}\label{pro1.1}
A quadratic polynomial with positive coefficients is stable.
\end{problem}
In general, the following partial converse holds:
\begin{problem}\label{pro1.2}
A polynomial of degree~$n$ with positive coefficients has no
roots in the sector
\begin{equation}\label{eq3}
|\arg z|\le \frac\pi n.
\end{equation}
\end{problem}
\noindent{\bf Hint:} Consider a broken line with $n{+}1$ segments
whose $k$th segment is parallel to the vector $a_kz^{n-k}$ 
($k=0,\dots,n$). If $\arg z$ is too small, this broken line cannot 
be closed.

\begin{problem}\label{pro1.3}
A polynomial
$p(z)=a_0z^n+a_1z^{n-1}+\dots+a_{n-1}z-a_n$
\
$(a_0,\dots,a_n>0)$ has exactly one root on the positive half-axis;
it is smaller than the absolute value of any other root of $p$.
\end{problem}

{\resume
The question how the signs of the coefficients affect the root 
distribution of a polynomial is quite interesting {\it per se\/} 
(see~\cite[part V, Chapter I]{1}), but leads away from the topic
of stable polynomials. On the other hand, the constructions of our 
next section turn out to be very fruitful.} 

\section{Nyquist-Mikhailov hodograph}\label{sec2}
Let $p(z)$ be a polynomial%
\footnote{For now, we do not need to assume that the coefficients
of the polynomial are real, although this is indeed the case in
most applications.}
of degree $n$.
In the complex plane $\mathbb{C}$, consider the curve%
\footnote{The normalizing factor~$i^{-n}$ is not of vital importance. It is
needed to simplify some formul{\ae} of the next sections.}
\begin{equation}\label{eq4}
\Gamma_p\equiv\left\{\,i^{-n}p(i\omega)\,\colon\,\omega\in\R\,\right\}.
\end{equation}
As the parameter~$\omega$ runs from $-\infty$ to $\infty$,
the curve is traversed in a certain direction. This
oriented curve is called the Nyquist-Mikhailov\footnote{Nyquist's name was added in translation [translators' remark].}  {\it hodograph}%
\footnote{The work of A. V. Mikhailov (1937) as well as the earlier work of
the American engineer H. Nyquist (1932) attracted attention to the geometrical
method that we describe here. Especially important
applications of this method were found in automatic control theory,
mostly thanks to papers of the Romanian mathematician V. M. Popov. However,
<<Mikhailov's hodograph>> was discovered by C. Hermite.}, or simply the
hodograph,  or the amplitude-phase characteristic of the polynomial $p$.

Assume that the polynomial $p$ has no roots on the imaginary axis. In
this case, $\Gamma_p$ does not go through zero and the function
\begin{equation}\label{eq5}
\phi_p(\omega)\equiv\Arg i^{-n}p(i\omega)=\im\Log i^{-n}p(i\omega)\quad
(\omega\in\R)
\end{equation}
is {\it continuous\/} at each point of the real axis.
Note that this function is defined up to an additive constant of the
form $2\pi k$ ($k\in\mathbb Z$), and its values do not have to lie
in the interval $[0,2\pi]$. In the sequel, we will be interested
in the increment%
\footnote{Along the way, we will also prove that the limits
$\phi_p(\pm\infty)$ exist.}
\begin{equation}\label{eq6}
\Delta_p\equiv\bigl.\phi_p\bigr|_{-\infty}^{+\infty},
\end{equation}
which is defined unambiguously.

\begin{lem}\label{lem2.1}
If\ \ $p(z)=z-\lambda$\ \ ($\re\lambda\ne0$),\ \ then\ \
$\Delta_p=-\pi\,\sgn\re\lambda$.
\end{lem}

\begin{proof}
The hodograph $\Gamma_p$ is a horizontal line traversed from left to right
that intersects the imaginary axis at the point $i\re\lambda$. Obviously,
 as $\omega$ runs from $-\infty$ to $+\infty$, the radius-vector of a
point on the hodograph makes a clockwise turn of magnitude $\pi$ if
$\re\lambda>0$  (counter-clockwise if $\re\lambda<0$).
\end{proof}

\begin{thm}[Hermite]\label{thm2.1}
If a polynomial $p$ has $n_{-}$ roots in the left half-plane and  $n_{+}$
roots in the right half-plane but no roots on the imaginary axis,
then
\begin{equation}\label{eq7}
\Delta_p=\pi(n_--n_+).
\end{equation}
\end{thm}
\begin{proof}
If $p=a_0\,p_1\,\cdots p_n$, where $p_k(z)=z-\lambda_k$, then
\begin{equation*}
\phi_p=\arg a_0+\phi_{p_1}+\dots+\phi_{p_n}\quad\text{and}\quad
\Delta_p=\Delta_{p_1}+\dots+\Delta_{p_n}.
\end{equation*}
From Lemma~\ref{lem2.1}, it follows that
$$
\Delta_p=-\pi\left(\sgn\re\lambda_1+\dots+\sgn\re\lambda_n\right)=
\pi(n_--n_+).
$$
\end{proof}

\begin{note}\label{not2.1}
Since there are no roots on the imaginary axis, we have $n_-+n_+=\deg p$.
Together with~\eqref{eq7}, this enables us to find both numbers $n_-$ and
$n_+$.
\end{note}

\begin{note}\label{not2.2}
The increment $\Delta_p$ achieves its maximal value, which equals~$\pi\deg p$,
for stable polynomials.
\end{note}

\begin{problem}\label{pro2.1}
If a polynomial $p$ is stable, prove that $\phi_p$ is monotone increasing
on $\R$.
\end{problem}

\begin{problem}\label{pro2.2}
Let $\gamma$ be a closed oriented curve on the Riemann sphere
(e.g., the imaginary axis is such a curve). Let us introduce the
<<generalized hodograph>>
$
\Gamma_p^\gamma\equiv
\left\{\,i^{-n}p(\omega)\,\colon\,\omega\in\gamma\,\right\}
$
and define the quantities $\phi_p^\gamma$ and $\Delta_p^\gamma$
analogously to~\eqref{eq5} and~\eqref{eq6}. Consider the unit circle
and the sector~\eqref{eq3} and formulate for them the analogue of
Theorem~\ref{thm2.1}.
\end{problem}

\begin{problem}\label{pro2.3}
Theorem~\ref{thm2.1} is related to the <<argument principle>> in
the theory of analytic functions. What is the value of the integral
$
\dfrac1{2\pi i}\int\limits_{-i\infty}^{+i\infty}\dfrac{p^\prime(z)}{p(z)}\,dz,
$
taken over the imaginary axis in the sense of the Cauchy principal value?
\end{problem}

\begin{problem}\label{pro2.4}
Investigate the hodograph of a rational function. Which rational functions
do you think should be called <<stable>>?
\end{problem}

\begin{problem}\label{probl2.5}
On the front page you see a stylized picture of the  hodograph of
the polynomial
$
p(z)=32z^6+12z^5+46z^4+21z^3+16z^2+7z+1
$
(the thickness of the curve decreases as the parameter $\omega$ increases).
Where are the roots of $p(z)$ located with respect to the imaginary axis?
Using {\sc Maple}, draw this curve on your own. Try to change one of
the coefficients of $p(z)$. What happens with the hodograph?
\end{problem}

{\resume
 Theorem~\ref{thm2.1} obtained by such simple tools can  already be applied
to count the number of roots of the polynomial $p$ to the left and to the right
of the imaginary axis. We can entrust the  drawing of hodographs to a computer and  determine the number of half-turns visually.

It is worth considering the following questions:
\newline\noindent
Why not use a computer to count all roots of a polynomial as well?
Which of the two problems will require more calculations?
For which polynomials will the pertinent calculations be hard and
their results unreliable? Which problems may arise for a developer and for
a user  of such a program?

Most likely, we will come to the conclusion that it is premature to
write a program, and the <<specific>> question of counting the number
$\Delta_p$ of half-turns  requires a mathematical rather than a
programming  solution. This is indeed the case.}

\section{Cauchy indices}\label{sec3}

The quantity $\Delta_p$ characterizes quite general topological properties
of the curve $\Gamma_p$. It turns out that for counting $\Delta_p$ it is
enough to know in which order a point moving along $\Gamma_p$
crosses the coordinate half-axes.

Let
\begin{equation}\label{eq8}
p(z)\equiv a_0z^n+a_1z^{n-1}+\dots+a_n\qquad(a_0\in\R,\ a_0>0).
\end{equation}
Let us consider the real polynomials
\begin{equation}\label{eq9}
\eqalign{
&f_0(\omega)\equiv
            +\re\left[i^{-n}p(i\omega)\right]=a_0\omega^n+\cdots,\cr
&f_1(\omega)\equiv
           -\im\left[i^{-n}p(i\omega)\right]=(\re a_1)\omega^{n-1}+\cdots}
\end{equation}
satisfying
\begin{equation}\label{eq10}
i^{-n}p(i\omega)=f_0(\omega)-if_1(\omega),\qquad\deg f_1<\deg f_0.
\end{equation}
If all coefficients of the polynomial~\eqref{eq8} are real,
then we have:
\begin{equation}\label{eq11}
\eqalign{
&f_0(\omega)\equiv a_0\omega^n-a_2\omega^{n-2}+a_4\omega^{n-4}-\cdots,\cr
&f_1(\omega)\equiv a_1\omega^{n-1}-a_3\omega^{n-3}+a_5\omega^{n-5}-\cdots  .
}
\end{equation}
Now it is time  to discuss the assumption of the previous section that
the polynomial $p$ has no roots on the imaginary axis. How can we check
this condition? From~\eqref{eq10} one can see that, for $\omega\in\R$,
\begin{equation}\label{eq12}
p(i\omega)=0\quad\Leftrightarrow\quad f_0(\omega)=f_1(\omega)=0\quad
                 \Leftrightarrow\quad\CMD(f_0,f_1)(\omega)=0.
\end{equation}
Thus, we need to use the Euclidean algorithm to find the greatest common
divisor $d\equiv\gcd(f_0,f_1)$. In the simplest case we get  $d=1$.
Otherwise, we need to find out whether the polynomial $d$ has real roots.

\begin{note}\label{not3.1} As we will see later, the idea of using
the Euclidean algorithm is extraordinarily fruitful. Now we simply
ran into it and risk to pass it by, not noticing that it is key to solving
the entire problem. Can we at this stage  {\it guess,\/} perhaps
only {\it feel,\/} the value of this accidental idea to develop it
afterwards?

Assume that the polynomial $p$ is continuously perturbed so that, at some
moment, one or several of its roots intersect the imaginary axis, thus
changing the values of $n_-$ and $n_+$, which we are interested in.
If we apply the Euclidean algorithm to the corresponding polynomials
$f_0$ and $f_1$, then its  {\it final result,\/} the greatest common
divisor $d$, will forget what happened. Most likely, $d$ was
equal to $1$ and will again become equal to $1$.
But what if the memory of those events will be preserved in the
{\it by-product\/} of the algorithm, which are usually thrown out
as useless? Later we will see that this guesswork will be
confirmed in its entirety.
\end{note}
Still assuming that the curve $\Gamma_p$ does not go through zero,
let us now consider how it intersects the imaginary axis. Let
$\omega_0<\omega_1<\dots<\omega_m$ be the values of the parameter
$\omega$ for which the intersections occur.  According to~\eqref{eq9},
$\omega_k$ are real roots of the polynomial $f_0$ of odd multiplicity
(so typically simple). Denote
\begin{equation}\label{eq13}
i_k\equiv
\lim_{\omega\to\omega_k}
\sgn\frac d{d\omega}\phi_p(\omega)\quad(k=0,1,\dots,m).
\end{equation}
In other words, $i_k=-1\quad(i_k=+1)$ if the radius-vector of a point
on the hodograph turns clockwise (counter-clockwise) when $\omega$
passes through the point $\omega_k$ (see the picture).

\begin{figure}[ht]\label{fig3.1}
\includegraphics[width=160mm, keepaspectratio=yes]{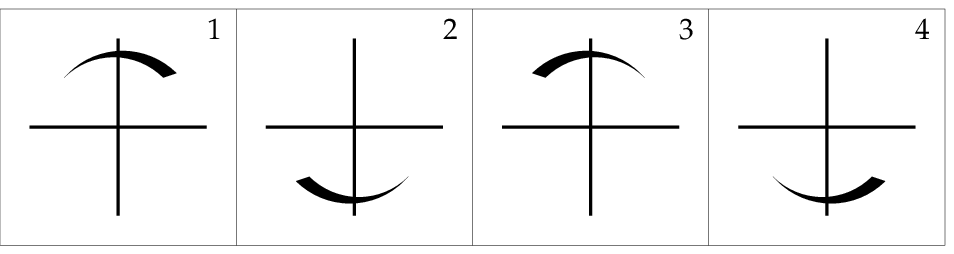}
\end{figure}

\begin{lem}\label{lem3.1}\quad
$\biggl.\phi_p\biggr|_{\omega_{k-1}}^{\omega_k}=
\dfrac\pi2\,\left(i_{k-1}+i_k\right)\quad(k=1,\dots,m)$
\end{lem}

\begin{proof}
For definiteness, suppose that $\omega_{k-1}$ corresponds
to an intersection of type~1 on the picture. Then the next
value~$\omega_{k}$ corresponds to an intersection of type~2 or~3.
In the former case, $\bigl.\phi_p\bigr|_{\omega_{k-1}}^{\omega_k}$ equals
$+\pi$, in the latter case,  zero, which agrees with the statement of
the lemma. The remaining three possibilities can be considered similarly.
\end{proof}

\begin{lem}\label{lem3.2}\quad
$\biggl.\phi_p\biggr|_{-\infty}^{\omega_0}=\dfrac\pi2\,i_0,\quad
\biggl.\phi_p\biggr|_{\omega_m}^{+\infty}=\dfrac\pi2\,i_m$
\end{lem}

\begin{proof}
Since $\deg f_0>\deg f_1$, we have
\begin{equation*}
\tan\phi_p(\omega)=-\frac{f_1(\omega)}{f_0(\omega)}\to 0\quad
                            (\omega\to\pm\infty).
\end{equation*}
Consequently, the limit directions of the radius-vector are horizontal.

For definiteness, suppose that the radius-vector approaches the direction
of the positive  real half-axis when $\omega\to-\infty$. Then between $-\infty$ and
$\omega_0$ there must be an intersection of type~2 or~3 (see the picture). In the former case
$\biggl.\phi_p\biggr|_{-\infty}^{\omega_0}$ equals
$+\dfrac\pi2$. In the latter case it is equal to~$-\dfrac\pi2$. This
agrees with the statement of the lemma. The remaining possibilities must
be considered similarly.
\end{proof}

\begin{lem}\label{lem3.3}
If $p$ has no roots on the imaginary axis, then
\begin{equation}\label{eq14}
\Delta_p=\pi\left(i_0+i_1+\dots+i_m\right)
\end{equation}
\end{lem}

\begin{proof}
The increment is additive:
\begin{equation*}
\Delta_p=\biggl. \phi_p\biggr|_{-\infty}^{\omega_0}+
         \biggl. \phi_p\biggr|_{\omega_0}^{\omega_1}+\dots+
         \biggl. \phi_p\biggr|_{\omega_{m-1}}^{\omega_m}+
         \biggl. \phi_p\biggr|_{\omega_m}^{+\infty}.
\end{equation*}
Therefore, the application of Lemma~\ref{lem3.1} and Lemma~\ref{lem3.2}
implies~\eqref{eq14}.

Strictly speaking, we should also consider the case $m=-1$, i.e., the case
when the hodograph does not intersect the imaginary axis at all.
Using the same reasoning as in the proof of Lemma~\ref{lem3.2}, we then will
see that~$\Delta_p=0$. This is indeed the proper way to understand the
formula~\eqref{eq14} in case $m=-1$.
\end{proof}

The quantities~$i_k$ and their sum occur in other applications. They have
specific names.

Let us consider the rational function
\begin{equation}\label{eq15}
R(\omega)\equiv\frac{f_1(\omega)}{f_0(\omega)}\qquad(\deg f_1<\deg f_0),
\end{equation}
where the numerator and the denominator are arbitrary real polynomials
and are not necessarily determined from~\eqref{eq9}.

Let $\omega_0<\omega_1<\dots<\omega_m$ be the real poles of~$R$ of odd order.
This means that $R(\omega)$ changes its sign as it <<goes through $\infty$>>
when $\omega$ goes through $\omega_k$.

The quantity
\begin{equation}\label{eq16}
\Ind\nolimits_{\omega_k}(R)\equiv
\begin{cases}
   &+1,\quad\text{if}\quad R(\omega_k-0)<0<R(\omega_k+0), \\
   &-1,\quad\text{if}\quad R(\omega_k-0)>0>R(\omega_k+0)
\end{cases}
\end{equation}
is called the {\it index\/} of the function~$R$ at its real pole $\omega_k$
of odd order.

The quantity
\begin{equation}\label{eq17}
\Ind\nolimits_a^b(R)\equiv\sum\limits_{k\,\colon\,a<\omega_k<b}
\Ind\nolimits_{\omega_k}(R)
\end{equation}
is called the {\it Cauchy index\/} of the function~$R$ on the interval~$(a,b)$.

A comparison with~\eqref{eq13} and with the picture now shows that
\begin{equation*}
i_k=\Ind\nolimits_{\omega_k}(R),\quad\text{where}\quad R\equiv\frac{f_1}{f_0}.
\end{equation*}

\begin{thm}\label{thm3.1}
Let the polynomial~\eqref{eq8} have $n_-$ roots in the left half-plane,
$n_+$ roots in the right half-plane, and no roots on the imaginary axis.
Then
\begin{equation*}
n_--n_+=\IndC\left(\frac{f_1}{f_0}\right),
\end{equation*}
where~$f_0$ and~$f_1$ are defined in~\eqref{eq9}.
\end{thm}

\begin{proof}
The theorem follows from Theorem~\ref{thm2.1} and from Lemma~\ref{lem3.3}.
\end{proof}

\begin{note}\label{not3.2}
Undoubtedly, the case when the polynomial is stable deserves special
consideration. We will devote to it a section of Chapter II.
\end{note}

\begin{note}\label{not3.3}
If we introduce the step function
\begin{equation*}
U(\omega)\equiv\sum\limits_{k\,\colon\,-\infty<\omega_k<\omega}
\Ind\nolimits_{\omega_k}(R),
\end{equation*}
then~\eqref{eq17} implies
\begin{equation}\label{eq18}
\Ind\nolimits_a^b(R)=\int\limits_a^b\,dU(\omega)=U(b-0)-U(a+0).
\end{equation}
If we could calculate values of the function~$U(\omega)$ without
calculating the roots $\omega_k$, we could conveniently apply
Theorem~\ref{thm3.1}.
\end{note}

\begin{problem}\label{pro3.1}
What is the Cauchy index of the logarithmic derivative
\begin{equation*}
R\equiv\frac d{d\omega}\ln f=\frac{f^\prime}f.
\end{equation*}
of a {\it real\/} polynomial $f$?
\end{problem}

\noindent{\bf Hint:}
It equals the number of {\it distinct\/} roots of the polynomial
$f$ in the interval $(a,b)$.

\begin{problem}\label{pro3.2}
The choice of the imaginary axis for computing values of $\Delta_p$
is most convenient. But we could take any other line that goes
through $0$, except for the real axis. Prove this and explain the
 problem with the real axis.
\end{problem}

\begin{problem}\label{pro3.3}
Let us draw a graph of the rational function $R$ on the {\it torus\/}
$\mathbb S^1\times\mathbb S^1$ instead of the plane $\R^1\times\R^1$.
This is useful since one of the points on the circle $\mathbb S^1$ must
correspond to the point at infinity on the axis $\R^1$. Verify that
such a <<graph>> is a closed curve on the torus. Which geometric
(more exactly, {\it topological\/}) meaning does $\IndC(R)$ acquire?
\end{problem}

{\resume
Theorem~\ref{thm3.1} is a better tool than Theorem~\ref{thm2.1}.
However, we seem to go around in circles:  trying  to avoid an
explicit calculation of the roots of the initial polynomial $p$,
we came to the necessity of calculating the roots of another polynomial
$f_0$! The next section will show that this is in fact unnecessary.}


\section{Sturm method}\label{sec4}

In the book~\cite{5}, one can find an elegant theorem of Sturm
on counting the number of real roots of a polynomial in a given interval.
In fact, this theorem --- or, more exactly, its method --- has a wider scope.

A finite sequence of polynomials $\{f_0,f_1,\dots,f_n\}$ is called
a {\it Sturm sequence\/} on the interval $(a,b)$ if
\begin{align}
&f_0(c)=0\quad(a<c<b)\quad\Rightarrow\quad f_1(c)\ne0 ,
\label{eq19}\\
&f_n(c)\ne0,\quad\forall c\in(a,b) ,
\label{eq20}\\
&f_k(c)=0\quad(0<k<n,\ a<c<b)\quad\Rightarrow\quad
 f_{k-1}(c)\,f_{k+1}(c)<0 .
\label{eq21}
\end{align}

Given a Sturm sequence, let us introduce the integer-valued function $V(x)$
defined to be the  number of sign changes in the sequence
\(
\{f_0(x),f_1(x),\dots,f_n(x)\}.
\)
The domain of this function is the interval $(a,b)$ from which the
roots of the polynomials in the sequence are excluded. At these points,
the function $V(x)$ can have a discontinuity of the first kind. However,
$V(x)$ does not have too many discontinuities:
\begin{lem}\label{lem4.1}
If $c\in(a,b)$ is not a zero of odd  multiplicity of the initial polynomial
$f_0$, then $V(c+0)=V(c-0)$.
\end{lem}

\begin{proof}
If $f_k(c)=0$, then $k<n$ in accordance with~\eqref{eq20}. Now let
 $k>0$; then~\eqref{eq21} implies that~$f_{k-1}(c)$
and~$f_{k+1}(c)$ are nonzero and have different signs.
Therefore, the subsequence $\{f_{k-1}(c),f_k(c),f_{k+1}(c)\}$
($0<|x-c|<\epsilon$) contains exactly one sign change regardless of
the sign of $f_k(x)$ for $x$ in a small neighborhood of $c$.

But if $c$ is a zero of the polynomial $f_0$ of even multiplicity,
then $f_0(x)$ does not change sign in the punctured neighbourhood
$0<|x-c|<\epsilon$. By \eqref{eq19}, this also applies to  $f_1(x)$.
Thus, the subsequence $\{f_0(x),f_1(x)\}$ has the same number of sign
changes ($0$ or $1$) to the left of the point $c$ as it does to
the right of $c$.
\end{proof}

\begin{lem}\label{lem4.2}
If $c\in(a,b)$ is a zero of $f_0$ of odd multiplicity,  then
\begin{equation}\label{eq22}
V(c+0)-V(c-0)=-\Ind\nolimits_c\left(\frac{f_1}{f_0}\right).
\end{equation}
\end{lem}

\begin{proof}
If the index at $c$ is equal to $+1$, then the function $f_1/f_0$ changes
its sign from $-$ to $+$ when $x$ goes through the point $c$. The subsequence
$\{f_0(x),\ f_1(x)\}$ thus loses a sign change and $V(x)$ decreases
by $1$. In case the index is negative, the opposite holds.
\end{proof}

\begin{thm}[Sturm]\label{thm4.1}%
\footnote{We foresaw the existence of such a formula --- see
Remark~\ref{not3.3}. It would resemble~\eqref{eq18} even more if we
defined $V(x)$ as the number of {\it sign retentions\/} rather
than the number  of sign changes.}
\begin{equation}\label{eq23}
\Ind\nolimits_a^b\left(\frac{f_1}{f_0}\right)=V(a+0)-V(b-0).
\end{equation}
\end{thm}

\begin{proof}
Use the two previous lemmata and the fact that the full increment
of a step function is the sum of its increments at points of
discontinuity.
\[
V(b-0)-V(a+0)=\sum\,[V(c+0)-V(c-0)],
\]
where summation is taken over all discontinuities of the function~$V$,
that is, over all zeros of $f_0$ of odd multiplicity that  lie in $(a,b)$.
\end{proof}

\begin{note}\label{not4.1}
Note that neither in the definition of Sturm sequences nor in the proof of
the Sturm theorem was it essential that $f_0,\dots,f_n$ are polynomials.
Only the fact that these functions are continuous and have finitely many
 roots was used\footnote{Thus, the Sturm method solves problems of {\it algebra\/} by means of {\it analysis.\/} As a result, it is disliked
by both algebraists and analysts.}.
\end{note}

Let us consider an important special case of Theorem~\ref{thm4.1}.

A Sturm sequence $\{f_0,f_1,\dots,f_n\}$ is called {\it regular\/} if
\[
\deg f_k=n-k\quad (k=0,1,\dots,n).
\]

\begin{thm}\label{thm4.2}
Let $h_k$ be the leading coefficient of a polynomial $f_k$ from
a regular Sturm sequence $\{f_0,f_1,\dots,f_n\}$. Then
\begin{equation}\label{eq24}
\IndC\left(\frac{f_1}{f_0}\right)=
n-2\,v(h_0,h_1,\dots,h_n),
\end{equation}
where $v(h_0,h_1,\dots,h_n)$ is the number of sign changes in the sequence
$\{h_0,h_1,\dots,h_n\}$.
\end{thm}

\begin{proof}
For large $|x|$, the sign of $f_k(x)$ coincides with the sign of its
leading term $h_k\,x^{n-k}$. Therefore,
\begin{equation*}
\eqalign{
&V(+\infty)=v(h_0,h_1,\dots,h_n),\cr
&V(-\infty)=v(h_n,-h_{n-1},\dots,(-1)^nh_0)=n-v(h_0,h_1,\dots,h_n)}
\end{equation*}
(the latter quantity is the number of {\it sign retentions,\/}
which, together with the  number of sign changes in the sequence
$\{h_0,\dots,h_n\}$, sums up to $n$).
\end{proof}

It remains to discuss how a Sturm sequence can be constructed
from an initial pair of polynomials~$\{f_0,f_1\}$.
For the case $\deg f_0>\deg f_1$, which we focus on,  a modified
Euclidean algorithm can be used. The modification
consists in changing the sign of the remainder at each step:
\begin{equation}\label{eq25}
f_{k-1}=d_kf_k-f_{k+1},\quad\deg f_{k+1}<\deg f_k.
\end{equation}
It does not change the original meaning of the algorithm since, at the last
step, we will still obtain $\gcd(f_0,f_1)$, but it guarantees that the
condition~\eqref{eq21} is satisfied. Regarding the conditions~\eqref{eq19}
and~\eqref{eq20}, for this construction they are equivalent (prove it).

\begin{note}\label{not4.2}
This is not the only possible method for constructing a Sturm sequence.
\end{note}

\begin{note}\label{not4.3}
<<Though the Sturm method is excellent in theory, it is not convenient
in practice  due to the enormous number of numerical coefficients
of various powers of $x$ in a sequence of Sturm functions when an
equation of high enough degree is given.>> (P.~L.~Chebyshev)
\end{note}

\begin{problem}\label{4.1}
How can we construct a Sturm sequence if~$\deg f_0\le \deg f_1$?
\end{problem}

\begin{problem}\label{4.2}
Use the result of Problem~\ref{pro3.1} and suggest an algorithm for counting
the number of {\it distinct\/} roots of a polynomial~$f$ on an interval
$(a,b)$. Use this algorithm to count the number of roots of the polynomial
$p$ on the imaginary axis.
\end{problem}

{\resume
We have solved the Routh-Hurwitz problem, and not only for real but
also for complex polynomials. We found out that this problem is
algorithmically equivalent to finding the greatest common divisor
of two polynomials and is therefore rather simple. It turns out
 that, when used to examine the stability of real polynomials,
it simplifies further due to the specific structure of the initial
polynomials~\eqref{eq11}. This will be taken up in the next section.}

\section{Routh scheme}\label{sec5}

Let polynomials $f_0$ and $f_1$ be defined as in \eqref{eq11}.
If $a_1\ne 0$, then the quotient and the remainder in \eqref{eq25} are
\begin{align}
&d_1(\omega)=c\,\omega\quad \left(c=\frac{a_0}{a_1}\right),
\label{eq26}\\
&f_2(\omega)=\left(a_2-c\,a_3\right)\,\omega^{n-2}-
             \left(a_4-c\,a_5\right)\,\omega^{n-4}+\cdots .
\label{eq27}
\end{align}

We see that, first of all, $f_2$ has the same structure as $f_0$ and $f_1$,
and if its leading coefficient is nonzero, then the same procedure can be
applied to the pair $\{f_1,f_2\}$. Secondly, the coefficients of $f_2$ occur
in the second row of the rectangular matrix
\begin{equation}\label{eq28}
\begin{pmatrix}
a_1\ a_3\ a_5\ \dots\\
a_0\ a_2\ a_4\ \dots
\end{pmatrix}
\end{equation}
after the Gaussian elimination of the entry $a_0$.

This is the basis for the computational  {\it Routh scheme.\/}
In textbooks and handbooks, the Routh scheme is usually described in
a form suitable for computations by hand. Thanks to the progress
of programming, mathematics now has new tools for recording its
algorithms. Let us use standard Pascal.

Let the coefficients of the polynomial~\eqref{eq8} be stored in the array
{\tt var h: array[0..n] of real;} The transition from $\{f_0,f_1\}$
to $\{f_1,f_2\}$ described in \eqref{eq26}-\eqref{eq27} corresponds
to the formal transition from the polynomial $p$ to the polynomial
\[
a_1\,z^{n-1}+\left(a_2-ca_3\right)\,z^{n-3}+a_3\,z^{n-3}+
             \left(a_4-ca_5\right)\,z^{n-4}+\dots,
\]
whereupon the process goes on provided that corresponding coefficients
are not zero. This algorithm is realized by the routine \Routh, which
returns the logical value {\tt true} after a normal completion; in this
case, it places the leading coefficients $h_0,h_1,\dots,h_n$ of the
Sturm sequence polynomials into the array {\tt h}, to which
Theorem 4.2 is then applied.

\begin{flushleft}{\tt
\spec{15}{\pmb{function} Routh(\pmb{var} h:array[0..n] of real):boolean;}\\
\spec{15}{\pmb{var} k,j:integer;}\\
\spec{15}{\phantom{\pmb{var} k,}c:real;}\\
\spec{15}{\pmb{begin}}\\
\spec{25}{k:=1;}\\
\spec{25}{\pmb{while} (k<n-1) and (h[k]<>0) \pmb{do}}\\
\spec{35}{\pmb{begin}}\\
\spec{45}{c:=h[k-1]/h[k];}\\
\spec{45}{k:=k+1;}\\
\spec{45}{j:=k;}\\
\spec{45}{\pmb{repeat}}\\
\spec{55}{h[j]:=h[j]-c*h[j+1];}\\
\spec{55}{j:=j+2}\\
\spec{45}{\pmb{until} j>=n}\\
\spec{35}{\pmb{end};}\\
\spec{25}{Routh:=(k=n-1)}\\
\spec{15}{\pmb{end} \{Routh\};}}
\end{flushleft}

After the $k$th step, <<the tail>>~{\tt[k..n]} of the array~{\tt h} is filled
with the alternating coefficients of the polynomials $\{f_k,f_{k+1}\}$, but
<<the head>> contains the leading coefficients of all preceding
polynomials of the sequence as some sort of <<useful trash>>.

We begin our analysis of the algorithm by criticizing it. If
{\tt Routh(h)=false}, then we will only know that the given polynomial
generates a nonregular sequence, but the question of where its roots are
located will remain open. This is not a dead end since the original
Sturm method can work with these nonregular situations too. We will not
go into details, since all is well in the case of interest to us:
\begin{thm}\label{thm5.1}
A polynomial is stable if and only if \/
{\tt Routh(h)=true} and $h_0$, $h_1,\dots$, $h_n$ are not zero
and of the same sign.
\end{thm}

\begin{proof}
By Theorems~\ref{thm3.1} and~\ref{thm4.1}, stability implies
 $V(-\infty)-V(+\infty)=n$. On the other hand, since $V(-\infty)\le n$
and $V(+\infty)\ge 0$, this must be the extreme case $V(-\infty)=n$,
$V(+\infty)=0$. The former equality shows that the length of a Sturm
sequence is maximal and equal to $n{+}1$; such a sequence is regular.
The latter equality implies that all $h_k$ are of the same sign.
Hereby the necessity is proved. Now, Theorems~\ref{thm4.2}
and~\ref{thm3.1} imply the sufficiency.
\end{proof}

Here comes a pleasant surprise of the algorithm%
\footnote{%
Why keep talking about sad things?

\global\let\R\originalR
<<\R{לאָמיר בעסער רעדן  פֿון עפּעס פֿרײלעכערס : װאָס הערט זיך עפּעס מכּח דער חלירע אין אָדעס? \ldots}

[Let us talk of something amusing.
What is the news about cholera in Odessa?]>>  (Sholom Aleichem)\global\let\R\newR}.
As \Routh \/  was developed, it was assumed that the polynomial has no roots
on the imaginary axis and that it generates a regular Sturm sequence
\[
f_k(z)=h_k\,z^{n-k}+\cdots\qquad(h_k\ne0,\quad k=0,1,\dots,n).
\]
But \Routh \/ {\it does not require division by the last two coefficients\/}
$h_n$ and $h_{n-1}$, so these coefficients, unlike the rest, can take
the value zero. Thus, the scope of \Routh \/ is wider than originally intended.
This property of the algorithm comes in handy. The point is that the loss
 of stability of a fixed point of the system~\eqref{eq1}%
\footnote{Under perturbation of the system~\eqref{eq1} and therefore of $A$ [translators' remark].} is accompanied
by effects determined by precisely how the eigenvalues of the matrix $A$ leave
the left half-plane%
\footnote{A new fixed point or a {\it limit cycle\/} can branch from
a fixed point that loses stability. {\it Bifurcation theory\/} studies
such phenomena.}.
The following two scenarios are most common:
\begin{itemize}
\item a simple real eigenvalue crosses the imaginary axis at the point~$0$
when it enters the right half-plane;
\item a pair of simple non-real eigenvalues crosses the imaginary axis at
points      $\pm i\omega$ ($\omega\ne0$).
\end{itemize}

The Routh algorithm can distinguish between these two variants:
\begin{thm}\label{thm5.2}
Let \/ {\tt Routh(h)=true}. Then,
\begin{itemize}
\item[a)] if $h_{n-1}\ne0$, $h_n\ne0$, then the polynomial has no roots on the imaginary
          axis, and
\begin{equation}\label{eq29}
n_+=v(h_0,\dots,h_n),\quad n_-=n-v(h_0,\dots,h_n);
\end{equation}
\item[b)] if $h_{n-1}\ne0$, $h_n=0$, then the polynomial has one simple root on the
          imaginary axis at the point $0$, and
\begin{equation}\label{eq30}
n_+=v(h_0,\dots,h_{n-1}),\quad n_-=n-1-v(h_0,\dots,h_{n-1});
\end{equation}
\item[c)] if $h_{n-1}=0$, $h_{n-2}h_n<0$, then the polynomial has no roots on the imaginary
          axis, and
\begin{equation}\label{eq31}
n_+=v(h_0,\dots,h_{n-2})+1,\quad n_-=n-v(h_0,\dots,h_{n-2})-1;
\end{equation}
\item[d)] if $h_{n-1}=0$, $h_{n-2}h_n>0$, then the polynomial has two simple roots
          on the imaginary axis at the points $\pm i\omega$ $(\omega\ne0)$, and
\begin{equation}\label{eq32}
n_+=v(h_0,\dots,h_{n-2}),\quad n_-=n-2-v(h_0,\dots,h_{n-2});
\end{equation}
\item[e)] if $h_{n-1}=0$, $h_n=0$, then the polynomial has one double root on the
          imaginary axis at the point~$0$, and
\begin{equation}\label{eq33}
n_+=v(h_0,\dots,h_{n-2}),\quad n_-=n-2-v(h_0,\dots,h_{n-2}).
\end{equation}
\end{itemize}
\end{thm}

\begin{proof}\footnote{We recommend that the reader give a proof on his/her own.}
The last three polynomials of the Sturm sequence are:
\[
f_{n-2}(x)=h_{n-2}x^2-h_n,\quad f_{n-1}(x)=h_{n-1}x,\quad f_n(x)=h_n .
\]
In case (a), the sequence is regular. In the remaining cases, one or both
polynomials $f_{n-1}$, $f_n$ are identically zero. Recall that we deal
with the Euclidean algorithm, therefore the construction of $f_k$ must
be stopped as soon as there is division without remainder in~\eqref{eq25};
then the last nonzero polynomial $f_k$ is the greatest common divisor
of the initial polynomials up to a numeric factor. We have:
\begin{equation*}
\eqalign{
&d(\omega)=\frac1{h_{n-1}}f_{n-1}(\omega)=\omega\qquad\text{in case (b)},\cr
&d(\omega)=\frac1{h_{n-2}}f_{n-2}(\omega)=\omega^2-\frac{h_n}{h_{n-2}}\quad
                                       \text{in cases (c), (d), (e)}.}
\end{equation*}

The statements of the theorem about the number and the location of roots
of $p(z)$ on the imaginary axis follow from~\eqref{eq12}. It suffices to
verify~\eqref{eq29}--\eqref{eq33}.

In the <<safe>> case (c), the last polynomial $f_{n-2}$ is not equal to zero on
the real axis, and the system $\{f_0,f_1,\dots,f_{n-2}\}$ is a Sturm sequence.
Therefore we have to use Theorem~\ref{thm4.1} instead of Theorem~\ref{thm4.2}.

Cases (b), (d) and (e) require correction of the initial polynomial $p(z)$.
It should be divided by $i^l\,d(-iz)$, where $l=\deg d$, in order to remove
its purely imaginary roots. At the same time, $f_0$ and $f_1$ should be divided
by their greatest common divisor $d$. Their leading coefficients $h_0$ and
$h_1$ do not change since we normalized $d$ beforehand. Therefore the corrected
polynomial requires no recalculations. We just have to take into account that
its degree has decreased and to reuse Theorem~\ref{thm4.2}.
\end{proof}

\begin{note}\label{not5.1}
A substantial difference between $h_n$ and $h_{n-1}$ is that
the last component of the array does not get processed.  It remains
equal to the last term $a_n$ of the polynomial $p$. On the contrary,
the penultimate term of the array is subjected to the largest number
of arithmetical operations.
\end{note}

\begin{problem}\label{pro1}
Run \Routh manually for polynomials of degree $3$ and $4$ and find
necessary and sufficient conditions for the stability of these polynomials.
\end{problem}

\begin{problem}\label{pro2} Program the routine
\begin{flushleft}
\spec{15}{\pmb{function} IsStable(h:array[0..n] of real):boolean;}
\end{flushleft}
that inputs the array~{\tt h} of coefficients of a real polynomial $p$ and
returns {\tt true} if the polynomial $p$ is stable or~{\tt false} otherwise.
Program also the routine {\tt IsStableStodola} which differs from
{\tt IsStable} by testing a polynomial first for the Stodola condition.
Is such an improvement of {\tt IsStable} useful?
\end{problem}

\begin{problem}\label{pro3} The following source code is written in {\sc APL}
(see \cite{12}):%
\begin{flushleft}
\spec{15}{$\pmb{\nabla}$ B$\pmb{\,\leftarrow}$ WhatIsIt A}\\
\spec{20}{LOOP: ((1=1$\pmb{\downarrow}\pmb{\rho}$A)%
$\,\pmb{\lor}\,$(A[2;1]$\pmb{\leqslant}$0)) /EXIT}\\
\spec{20}{A$\,\pmb{\leftarrow}$%
(2 2$\pmb{\,\rho\,}$0 1 1 -A[1;1]$\pmb{\div}$A[2;1])%
$\ \pmb{+}\pmb{\centerdot}\pmb{\times}\ $(0 1$\,\pmb{\downarrow}\,$A)}\\
\spec{24}{$\pmb{\rightarrow}\,$LOOP}\\
\spec{20}{EXIT: B$\,\pmb{\leftarrow}$(A[2;1]>0)}\\
\spec{15}{$\pmb{\nabla}$}
\end{flushleft}
The input parameter {\tt A} is a matrix of the form~\eqref{eq28}. Which
algorithm does the function {\tt WhatIsIt} realize? How can it be improved?
 Give a comparative analysis of Pascal versus {\sc APL} as algorithmic
languages.
\end{problem}

{\resume
So we now have at our disposal a simple (possibly unimprovable) algorithm
for testing the stability of a real polynomial. We just have no
explicit formul{\ae}  that express the output of this algorithm.
For comparison, the Gauss algorithm and  the Cramer formul{\ae}
 are complementary in the theory of linear systems: the former describes
how to find a solution, the latter how this solution looks.}

\section{Hurwitz theorem}\label{sec6}

Given a  polynomial~\eqref{eq8} with real coefficients, consider a
corresponding $n{\times}n$ matrix of the following structure:
\begin{equation}\label{eq34}
\mathcal{H}_p\equiv
\begin{pmatrix}
a_1   &a_3   &a_5   &a_7   &\cdots\\
a_0   &a_2   &a_4   &a_6   &\cdots\\
0     &a_1   &a_3   &a_5   &\cdots\\
0     &a_0   &a_2   &a_4   &\cdots\\
0     &0     &a_1   &a_3   &\cdots\\
\vdots & \vdots & \vdots & \vdots & \ddots
\end{pmatrix}
\end{equation}
(the coefficients $a_0,\dots,a_n$ are not enough to fill the rows, but we
set $a_{n+1}=a_{n+2}=\dots=0$).
The matrix $\mathcal{H}_p$ is called the {\it Hurwitz matrix\/} of the
polynomial $p$. Let us denote by $\eta_k$ the leading principal minor
of this matrix formed from the first $k$ rows and columns. It is easy
to check that there is only one (diagonal) nonzero term in the last
column of $\mathcal{H}_p$. It equals $a_n$. Therefore
\begin{equation}\label{eq35}
\eta_n=\eta_{n-1}a_n.
\end{equation}

\begin{lem}\label{lem6.1}
In the regular case,
\begin{equation}\label{eq36}
h_1=\eta_1,\ h_2=\frac{\eta_2}{\eta_1},\ \dots,\
h_n=\frac{\eta_n}{\eta_{n-1}}.
\end{equation}
\end{lem}

\begin{proof}
First note that the matrix $\mathcal{H}_p$ consists of blocks  of
type~\eqref{eq28}, which are in turn made up from the coefficients of
the polynomials $f_0$ and $f_1$. Let us reduce the matrix
$\mathcal{H}_p$ to upper triangular form using Gaussian elimination
without pivoting. By \eqref{eq26}-\eqref{eq27}, the elimination
of the term $a_0$ from the second, fourth etc.\ rows leaves these
rows  filled with the  coefficients of $f_2$, the next polynomial
of the sequence. The first row and column are no longer needed.
Crossing them out (temporarily), we obtain an  $(n-1)\times(n-1)$
matrix of the same <<Hurwitz>> structure formed from the coefficients
of the polynomials $f_1$ and $f_2$. Repeating the same procedure,
 we will arrive at a  triangular matrix
\begin{equation*}
\begin{pmatrix}
h_1   & * & * &\cdots & * \\
0     &h_2    & * &\cdots & *\\
0     &0      &h_3    &\cdots & * \\
\vdots&\vdots &\vdots &\ddots &\vdots\\
0     &0      &0      &\cdots &h_n
\end{pmatrix}
\end{equation*}
whose $k$th row consists of the coefficients of the polynomial $f_k$,
starting with the leading coefficient $h_k$ on the main diagonal.

Each elementary transformation that we applied consisted in subtracting
the {\it preceding\/} row (multiplied by a suitable number) from a given
row. Such  transformations preserve not only the determinant
$|\mathcal{H}_p|=\eta_n$  but also all leading principal minors
$\eta_k$. Consequently,
\begin{equation}\label{eq37}
\eta_k=h_1h_2\cdots h_k\qquad(k=1,2,\dots n),
\end{equation}
which implies \eqref{eq36}.
\end{proof}

\begin{thm}[Hurwitz]\label{thm6.1} A polynomial
\begin{equation}\label{eq38}
p(z)=a_0z^n+a_1z^{n-1}+\dots+a_n\qquad(a_0,a_1,\dots,a_n\in\R;\ a_0>0)
\end{equation}
is stable if and only if all leading principal minors of its Hurwitz matrix
$\mathcal{H}_p$ are positive:
\begin{equation}\label{eq39}
\eta_k>0\qquad(k=1,2,\dots,n).
\end{equation}
\end{thm}

\begin{proof}
The assertion of the theorem follows immediately from the preceding Lemma
and Theorems~\ref{thm5.1} and \ref{thm5.2}.
\end{proof}

\begin{note}\label{not6.1}
Let us assume that the polynomial~$p$ depends continuously on one or several
parameters and is initially stable. Lemma 6.1 and Theorem 5.2 imply
that its  stability is maintained for as long as $|\mathcal{H}_p|=\eta_n$
is not zero%
\footnote{Unlike the fish, the Hurwitz conditions <<rot from the tail>>.}.
At the same time, owing to~\eqref{eq34} and to Theorem 5.2 again,
a sign change in $a_n$ (for $\eta_{n-1}\ne 0$) corresponds to the
transition of a simple real root from the left to the right half-plane.
A sign change in $\eta_{n-1}$ (for $a_n\ne0$, $\eta_{n-2}\ne 0$)
corresponds to the transition of a pair of simple complex conjugate
roots from the left to the right half-plane.
\end{note}

\begin{note}\label{not6.2}
The observation we just made appears useful: in order to <<catch>>
the moment when a continuously varying polynomial loses stability,
it is sufficient to check only the sign of the last minor $\eta_n$.
This is, however, no reason to celebrate:  Theorem \ref{thm5.1} implies
that, while the polynomial remains stable, the determinant $|\mathcal{H}_p|$
is best computed using the compact Gauss scheme, which  yields all the
leading principal minors $\eta_k$ as a by-product anyway.
In general, it makes no sense to apply the Hurwitz theorem for
computations -- to this end one should use the Routh scheme.
\end{note}

\begin{note}\label{not6.3}
In 1914,  A. Li\'enard and M. Chipart proposed a different criterion of
polynomial stability. They established that a polynomial \eqref{eq38}
of degree $n$ {\it with positive coefficients\/}  is stable if and only
if the following conditions are satisfied:
\begin{equation}\label{eq40}
\eqalign{
&\eta_2>0,\eta_4>0,\dots,\eta_{n-1}>0,\qquad\text{if}\ n\ \text{is odd},\cr
&\eta_1>0,\eta_3>0,\dots,\eta_{n-1}>0\qquad\text{if}\ n\ \text{is even}.}
\end{equation}

The Li\'enard-Chipart conditions \eqref{eq40} look simpler than the Hurwitz conditions
\eqref{eq39} since they contain half that many determinantal inequalities.
Although the simplicity is misleading from the computational point of view
(for reasons given in Remark \ref{not6.2}), conditions \eqref{eq40} may be
more useful for formal derivations, and the equivalence of \eqref{eq39}
and \eqref{eq40} is very interesting from the theoretical point of view.
\end{note}

\begin{problem}\label{pro5.1}
Assuming the minors $\eta_1,\dots,\eta_n$ are all nonzero, prove that
the number of zeros of the polynomial $p$ in the right (left) half-plane
is equal to the number of sign variations (retentions) in the sequence
\(
\left\{\,a_0,\ \eta_1,\ \dfrac{\eta_2}{\eta_1},\ \dots,\
                        \dfrac{\eta_n}{\eta_{n-1}}\,\right\}.
\)
\end{problem}

\begin{problem}\label{pro5.2}
Let all $\eta_k$ be nonnegative. Does this imply that
all roots of $p$ lie in the {\it closed\/} left half-plane?
\end{problem}

\noindent{\bf Hint:}
No, it does not. Give a counterexample.

\begin{problem}\label{pro5.3} Prove that the statement of the Hurwitz theorem
remains valid if $a_0,a_1,a_2,\dots$ in~\eqref{eq34} are the coefficients of
the polynomial
\begin{equation}\label{eq41}
p(z)\equiv a_0+a_1z+a_2z^2+\cdots\qquad(a_0,a_1,\ldots\in\R;\ a_0>0).
\end{equation}
\end{problem}

\noindent{\bf Hint:} if $p(0)\neq0$, then the polynomials $p(z)$ and
$q(z)\equiv z^n\,p(z^{-1})$ ($n=\deg p$) are simultaneously stable or
unstable.

\begin{problem}\label{pro5.4}
If $p(z)$ is an analytic function represented by the power series~\eqref{eq41}, then
we can formally construct its infinite Hurwitz matrix and require the positivity
of its leading principal minors. Give an example showing that the Hurwitz theorem
does not generalize to analytic functions. In this connection,
see \cite{2,6,7}.
\end{problem}

{\resume
Our goals are achieved, and we are done with the basic material.
In Chapter II we will consider many variations on the topic
of stable polynomials. That optional material can be used for
seminars as well as for independent study.}

\section*{Appendix to Chapter I}

Here we consider the stability problem of the Lorentz system
posed in the Introduction%
\footnote{We remind the reader that the parameters $b$, $r$, $\sigma$
are assumed to be positive [translators' remark].}:
\begin{equation}\label{lorentz}
    \dot{X}=\sigma Y-\sigma X,\quad
    \dot{Y}=rX-Y-XZ,\quad
    \dot{Z}=XY-bZ.
\end{equation}

\noindent {\bf 1. Determining fixed points.} A fixed point of the system~(\ref{lorentz})
is a point $(X_k,Y_k,Z_k)\in\R^3$ that annihilates the right hand sides in~(\ref{lorentz}), i.e., a solution to the following system of algebraic equations:
  \begin{equation}\label{l43}
    \sigma Y_k-\sigma X_k=0,\quad
    rX_k-Y_k-X_kZ_k=0,\quad
    X_kY_k-bZ_k=0.
\end{equation}
One of the solutions to~(\ref{l43}) is easy to spot: it is the zero fixed point
  \begin{equation}\label{l44}
   (X_0,Y_0,Z_0)=(0,0,0).
\end{equation}
In addition, the system~(\ref{l43}) has two more solutions
 \begin{equation}\label{l45}
   (X_{1,2},Y_{1,2},Z_{1,2})=(\pm\sqrt{b(r-1)}, \pm\sqrt{b(r-1)}, r-1).
\end{equation}

\begin{note} The zero fixed point exists for all positive values of the parameters
$\sigma$, $b$, $r$. A pair of nonzero fixed points~(\ref{l45}) bifurcates from it when $r>1$. \end{note}  \vskip 0.3cm

\noindent {\bf 2. Linearization.}   In order to linearize the system $\dot{x}=f(x)$ in the neighborhood of a fixed point $x_k$ (where $f(x_k)=0$), one has to replace the
function $f(x)$ by its linear form $A_k(x-x_k)$ where $A_k=f'(x_k)$ is the Jacobian
matrix consisting of the partial derivatives of the function $f(x)$ taken at the point
$x_k$. In our case,
\[  A_k = \left(  \begin{array}{ccc}
-\sigma & \sigma & 0 \\ r-Z_k & -1 & -X_k \\ Y_k & X_k & -b
\end{array}
\right)   \qquad  (k=0,1,2).
\]
\vskip 0.3cm

\noindent {\bf 3. Computing the characteristic polynomial.}  We have
\[
 p_k(\lambda)  =  | \lambda I - A_k |  =  (\lambda + \sigma) (\lambda + 1) (\lambda + b)
 +\sigma X_k^2+ \sigma (Z_k-r) (\lambda + b) + X_k^2 (\lambda + \sigma)
\]
since $X_k=Y_k$ for all three fixed points.

For the zero fixed point~(\ref{l44}), we obtain
\[
p_0(\lambda)  =  (\lambda + b) [\lambda^2 + (\sigma+1) \lambda +
\sigma (1-r)];
\]
for the fixed points~(\ref{l45}), we get
\begin{eqnarray*}
p_{1,2}(\lambda)  = a_0 \lambda^3 + a_1 \lambda^2 + a_2 \lambda  + a_3, \qquad
\text{where} \\
a_0=1, \;  \; a_1 =\sigma + b + 1 ,  \;  \; a_2 = b (\sigma+r), \; \;
a_3 = 2 \sigma b (r-1).
\end{eqnarray*} \vskip 0.3cm

\noindent {\bf 4. Stability of fixed points.} The Hurwitz theorem is in fact not needed
to analyze the stability of the polynomial $p_0$ -- the elementary Vieta's Theorem
is enough. This polynomial is stable for $0<r<1$.

\begin{note} At $r=1$, one of the roots of the polynomial $p_0$ crosses
the imaginary axis and enters the right half-plane. The zero fixed point
loses its stability, and exactly then the fixed points~(\ref{l45}) bifurcate
from it. \end{note}

For the polynomials $p_{1,2}$ of degree $3$, the Hurwitz matrix has the form
\[  \left(  \begin{array}{ccc}
a_1 & a_3 & 0 \\ a_0 & a_2 & 0 \\ 0 & a_1 & a_3
\end{array}
\right),
\]
and the conditions of Theorem~\ref{thm6.1} say
\[  a_1>0, \qquad
\left|  \begin{array}{cc}
a_1 & a_3 \\ a_0 & a_2 \end{array}
\right|>0, \qquad
 \left|  \begin{array}{cc}
a_1 & a_3  \\ a_0 & a_2
\end{array}
\right| a_3 >0.
\]
As $a_1>0$, $a_3>0$ for $r>1$,  there remains just one inequality
\[ a_1 a_2 - a_0 a_3 = b (\sigma + b + 1) (\sigma + r) - 2 b \sigma (r-1) > 0,   \]
which can be viewed as an answer. However, it is better to represent the answer as
follows:
\[  1< r < r_*  \; \; \text{where}  \; \; r_* = \begin{cases}\sigma { \sigma + b + 3
\over \sigma - b - 1  }  & \sigma > b + 1, \\ +\infty & \text{otherwise}. \end{cases}  \]
\begin{note} The fixed points that bifurcate from $(0,0,0)$ are initially stable.
When the parameter $r$ reaches its critical value $r_*$, they lose their stability
since, by Theorem~\ref{thm5.2} (also see Remark~\ref{not6.1}), a pair of nonreal
roots of the polynomial $p_{1,2}$ crosses the imaginary axis. One can in fact show
that, at $r=r_*$, {\it limit cycles\/} branch from the fixed points that lose
stability. Further discussion of these difficult and interesting questions is beyond
the scope of our notes; these questions are subject of bifurcation theory.
\end{note}

\chapter{Extensions}
\vskip 3cm
\begin{figure}[ht]
\includegraphics[width=160mm, keepaspectratio=yes]{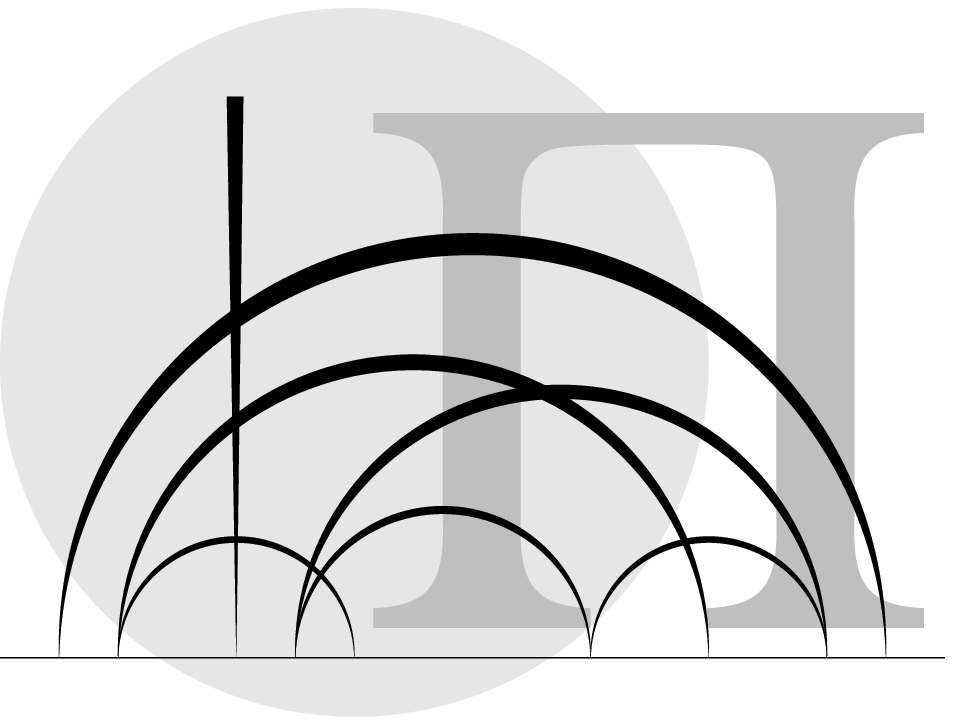}
\end{figure}
\begin{quotation}{\it

"Voil\`a le sujet simplifi\'e, argumentum omni denudatum ornamento. Je ferais avec cela, continua le j\'esuite, deux volumes de la taille de celui-ci."

Et, dans son enthousiasme, il frappait sur le saint Chrysostome in-folio qui faisait plier la table sous son poids.

D'Artagnan fr\'emit.}

\end{quotation}

\hfill
Alexander Dumas.  {\sc Les Trois Mousquetaires.} 

\hfill Chapitre XXVI. La th\`ese d'Aramis.

\vskip 2cm

\setcounter{section}{6}

\section{Stieltjes Fractions}\label{sec7}

The definition of the Cauchy index $\IndC(R)$ given in \S\,\ref{sec3}
presupposed that
the rational function $R$ vanishes at infinity. Let us try to generalize this
notion to arbitrary rational functions. A useful idea is contained in
Problem \ref{pro3.3}: the point at infinity $\infty$ should be <<made equal>> to
the other points, and it makes sense to consider the function $R$ as
a map of the {\it projective line\/}%
\footnote{For sure, the point is not just to add a new element to the
{\it set\/}  $\mathbb{R}^1$; we need to incorporate that element into
the algebraic and topological structures that exist in $\mathbb{R}^1$.
More about projective spaces $\PR^n$ can be found in \cite{13, 14}.}
$\PR^1\equiv\mathbb{R}^1\cup\{\infty\}$ into itself. In particular, if
$R=f_1/f_0$ and  $\deg f_1>\deg f_0$, then we will view $R$ as having
a pole at $\infty$ of order  $\deg f_1-\deg f_0$. If that pole is of
odd order, then, analogously to \eqref{eq16} (\S\,\ref{sec3}), 
we let
\begin{equation}\label{eq46}
\IndI(R)=
\begin{cases}
&+1,\quad\text{if}\quad R(+\infty)<0<R(-\infty)\\
&-1,\quad\text{if}\quad R(+\infty)>0>R(-\infty),
\end{cases}
\end{equation}
in all other cases, let $\IndI(R)=0$. This summand should be added to
\eqref{eq17} (\S\,\ref{sec3}):
\begin{equation}\label{eq47}
\IndPR(R)\equiv\IndC(R)+\IndI(R) .
\end{equation}
\begin{note}\label{not71}
Incidentally, \eqref{eq46} implies that rays $(C,+\infty)$ should
be viewed as {\it left\/} half-neighborhoods of the point $\infty$,
and rays $(-\infty,-C)$ as its {\it right\/} neighborhoods. The projective
line is in one-to-one correspondence with the circle $\mathbb S^1$,
as illustrated below:
\end{note}

\begin{figure}[ht]\label{fig71}
\includegraphics[width=160mm, keepaspectratio=yes]{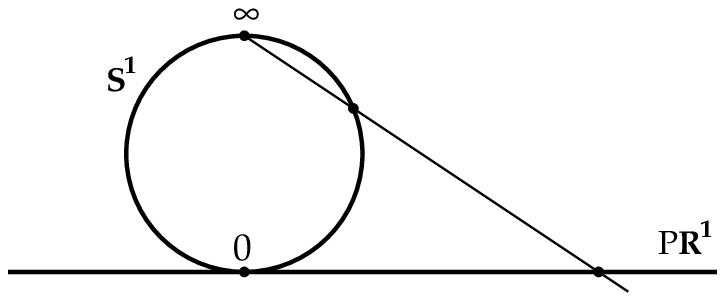}
\end{figure}

\begin{note}\label{not72}
Let $d(\omega)=c\omega^\nu+\cdots$ ($\nu=\deg d$) be a polynomial. Then
\begin{equation}\label{eq48}
\IndI(d)=
\begin{cases}
-\sgn c,\quad&\text{if}\quad\nu\quad\text{is odd,} \\
 0,     \quad&\text{if}\quad\nu\quad\text{is even.} \\
\end{cases}
\end{equation}
\end{note}

Now note how the index $\IndPR(R)$ changes when various <<projective>>
maps act on $R$.

\begin{lem}\label{lem71}
If $d\in\mathbb{R}$ is a constant, then $\IndPR(d+R)=\IndPR(R)$.
\end{lem}
\begin{proof}
The addition of a constant does not change the behaviour of a function near its poles.
\end{proof}
\begin{lem}\label{lem72}
If $d$ is a polynomial and $R(\infty)=0$, then $\IndPR(d+R)=\IndC(R)+\IndI(d)$.
\end{lem}
\begin{proof}
See \eqref{eq46}--\eqref{eq47}.
\end{proof}
\begin{lem}\label{lem73}
$\IndPR\left(-\tfrac1R\right)=\IndPR(R)$.
\end{lem}
\begin{proof} The projective line $\PR^1$ is divided by its {\it two\/}
points $0$ and $\infty$ into its positive $\mathbb{R}^+$ and negative
 $\mathbb{R}^-$ rays. Let a variable $\omega$ traverse $\PR^1$ and return
to its starting point. Clearly, the number of crossings from $\mathbb{R}^-$
to $\mathbb{R}^+$ must equal the number of reverse crossings. The crossings
through $\infty$ occur at the poles of the function $R$,
and they are accounted for in the sum \eqref{eq47} with the
appropriate sign.  The crossings through $0$ occur at the zeros of $R$,
i.e., at the poles of  $\tfrac1R$, and they are accounted for in the
analogous formula for $\IndPR\left(\tfrac1R\right)$. As a result,
$\IndPR\left(\tfrac1R\right)+\IndPR(R)=0$.
\end{proof}
\begin{note}\label{not73}
The transformations $\omega\mapsto\omega+d$ are called {\it shifts\/}
of the projective line $\PR^1$. The point $-\tfrac1\omega$ is referred
to as the {\it polar\/} of $\omega\in\PR^1$. If the diameter of the circle
on Fig.\,2 is equal to 1, then the polar $B$ of the point $A$ is constructed
by drawing the perpendicular $BC$ to $AC$. By a theorem from elementary
geometry,  $OA\cdot OB=OC^2$. Another theorem (on a subscribed and a
central angle) implies that $A^\prime B^\prime$ is a diameter.
\end{note}

\begin{figure}[ht]\label{fig72}
\includegraphics[width=160mm, keepaspectratio=yes]{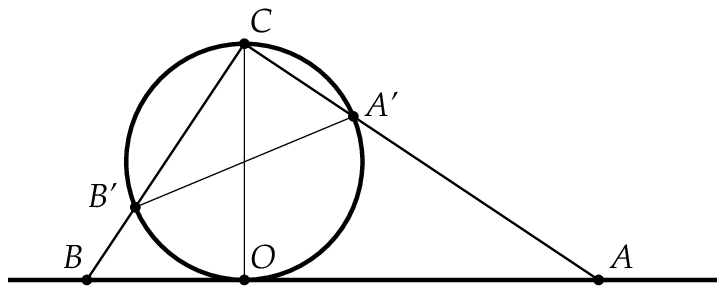}
\end{figure}

Alternating between shifts $\omega\mapsto\omega+d$ and
the polar transformation $\omega\mapsto-\tfrac1\omega$, we obtain a
{\it continued fraction\/}. The right-hand sides of the following
formula are examples of continued fractions:
\begin{lem}\label{lem74}
If\/ $\alpha\delta-\beta\gamma=1$, then
\begin{equation*}
\dfrac{\alpha\omega+\beta}{\gamma\omega+\delta}=
\begin{cases}
 \dfrac1\gamma-
 \dfrac1{\gamma-
 \dfrac1{\dfrac1\gamma-
 \dfrac1{\alpha\gamma-
 \dfrac1{\dfrac\delta\gamma+\omega}}}}\qquad&(\gamma\ne0), \\
 (\alpha\beta+\alpha)-
 \dfrac1{\dfrac1\alpha-
 \dfrac1{\alpha-
 \dfrac1\omega}}\qquad&(\gamma=0).
\end{cases}
\end{equation*}
\end{lem}
\begin{lem}\label{lem75}
If\/ $\alpha\delta-\beta\gamma=1$, then
\begin{equation}\label{eq49}
\frac{\alpha\omega+\beta}{\gamma\omega+\delta}=
\begin{cases}
\cfrac1{\gamma-
 \cfrac1{\gamma-
  \cfrac1{
   \cfrac1{\gamma-
    \dfrac1{\alpha\gamma-
     \cfrac1{\dfrac\delta\gamma+\omega
}}}}}}\qquad&(\gamma\ne0),\\
(\alpha\beta+\alpha)-
\cfrac1{
 \cfrac1{\alpha-
  \cfrac1{\alpha-
   \cfrac1{\omega
}}}}\qquad&(\gamma=0).
\end{cases}
\end{equation}
\end{lem}
\begin{proof}
By straightforward calculation.
\end{proof}

The fractional-linear map
\begin{equation*}
\Phi\ \colon\ \omega\mapsto\frac{\alpha\omega+\beta}
                                {\gamma\omega+\delta}\quad
(\alpha,\beta,\gamma,\delta\in\mathbb{R};\ \alpha\delta-\beta\gamma=1)
\end{equation*}
performs a {\it projective transformation\/} $\PR^1\to\PR^1$.
Such transformations form the {\it group\/}%
\footnote{One can find information about this group --
and other things -- in  \cite[\S 5]{15}.}
denoted by $\SL$.
\begin{thm}\label{thm71}
The Cauchy index is invariant under the action  of the group
$\SL$ on rational functions:
\begin{equation*}
\IndPR(\Phi\circ R)=\IndPR(R),\quad\forall\Phi\in\SL.
\end{equation*}
\end{thm}
\begin{proof}
Apply Lemmata~\ref{lem74}, \ref{lem71}, \ref{lem73}.
\end{proof}

Now let us explore connections between projective geometry, continued
fractions, and the Euclidean algorithm. For definiteness, suppose that
$R(\infty)=0$, so that $R=f_1/f_0$,\ $\deg f_1<\deg f_0$ in Sections~3 and~4.
Run the modified Euclidean algorithm, which we used in \S\,\ref{sec4}  to construct Sturm sequences:
\begin{equation}\label{eq50}
f_{k-1}=d_k f_k-f_{k+1}\quad
(k=1,\dots,m;\ \deg f_{k+1}<\deg f_k;\ f_{m+1}=0).
\end{equation}
If
\begin{equation}\label{eq51}
R_k\equiv\frac{f_{k+1}}{f_k}\quad
(k=0,1,\dots,m;\ R_0=R,\ R_m=0),
\end{equation}
then \eqref{eq50} yields the recurrence relation
\begin{equation}\label{eq52}
-R_{k-1}=-\frac1{d_k-R_k},
\end{equation}
and we obtain the following expansion of the rational function $-R$
into a continued fraction:
\begin{equation}\label{eq53}
-R=-\cfrac1{d_1-\cfrac1{d_2-\cfrac1{\ddots-\cfrac1{d_m}}}},
\end{equation}
A functional continued fraction of type~\eqref{eq53}, where
 $d_1,\dots,d_m$ are polynomials, is called a Stieltjes continued
fraction (see~\cite{1, 16a, 16b, 16c, 16d}).
\begin{thm}\label{thm72}
If a rational function $R$ is represented by a continued
fraction~\eqref{eq53}, then
\begin{equation}\label{eq54}
\IndPR(R)=-\sum\limits_{k=1}^m \IndI(d_k).
\end{equation}
\end{thm}
\begin{proof}
Apply Lemmata~\ref{lem72} and \ref{lem73} inductively, using the
relation~\eqref{eq52}.
\end{proof}
\begin{note}\label{not74}
Theorem~\ref{thm72} shows a way to compute Cauchy indices,
which parallels  Schur's method. It uses the same (Euclidean)
algorithm, but its validity is established by different reasoning.
\end{note}

Now consider the extreme case when the function $R$ is generated by
a stable polynomial, as described in  \S\,\ref{sec3}. 
It is time to give such functions a name. We will call a rational function $R=f_1/f_0$
{\it proper\/} if  $\deg f_1<\deg f_0$ and $\IndPR(R)=\deg f_0$;
the class of all proper functions will be denoted by $\mathcal R$.
\begin{thm}\label{thm73}
$R\in\mathcal R$ if and only if
\begin{equation}\label{eq55}
-R(\omega)=-
\cfrac1{\alpha_1\omega+\beta_1-
 \cfrac1{\alpha_2\omega+\beta_2-
  \cfrac1{\ddots-
   \cfrac1{\alpha_n\omega+\beta_n
}}}},
\end{equation}
where $\beta_1,\dots,\beta_n\in\mathbb{R},\ \alpha_1,\dots,\alpha_n>0$.
\end{thm}
\begin{proof}
Formula \eqref{eq50} implies that $n=\sum\nolimits_{k=1}^m\deg d_k$, whereas
\eqref{eq48} implies the inequality $-\IndI(d)\leq\deg d$, where equality
occurs if and only if  $d(\omega)=\alpha\omega+\beta$ ($\beta\in\mathbb{R},
\alpha>0$). A comparison with \eqref{eq54} now shows that the polynomials
$d_1,\dots,d_m$ in \eqref{eq53} must satisfy exactly these conditions,
and their number must be exactly $n$.
\end{proof}
\begin{note}\label{not75}
The transformations
$z\mapsto\alpha z+\beta$ ($\beta\in\mathbb{R},\ \alpha>0$) and
$z\mapsto-\dfrac1z$ map the complex {\it upper\/} half-plane
$\{z\in\mathbb{C}\ \colon \im z>0\}$ into itself. Therefore,
this property is inherited by the function  $-R$ if it admits
an expansion of type \eqref{eq55}. This theme will be taken up again
in \S\,\ref{sec8}, but our next problem already has 
a hint of a variation on this idea:
\end{note}
\begin{problem}\label{pro71}
Let
\begin{equation}\label{eq56}
R(\omega)=\dfrac1{\alpha_1\omega+\beta_1+
          \dfrac1{\alpha_2\omega+\beta_2+
          \dfrac1{\ddots+
          \dfrac1{\alpha_n\omega+\beta_n}}}},
\end{equation}
where $\alpha_k,\beta_k>0$ ($k=1,2,\dots,n$).
Then the function $R$ maps the {\it right\/} half-plane
$\{z\in\mathbb{C}\ \colon\ \re z>0\}$ into itself;
all its poles and zeros must lie in the {\it left\/} half-plane.
\end{problem}
\begin{problem}\label{pro72}
Consider a tridiagonal matrix
\begin{equation*}
A=
\begin{pmatrix}
a_1   &c_2   &0     &\dots  &0       &0      \\
b_2   &a_2   &c_3   &\dots  &0       &0      \\
0     &b_3   &a_3   &\dots  &0       &0      \\
\vdots & \vdots &  \vdots & \ddots &  \vdots & \vdots  \\
0     &0     &0     &\dots  &a_{n-1} &c_n    \\
0     &0     &0     &\dots  &b_n     &a_n
\end{pmatrix} .
\end{equation*}
Prove that
\begin{itemize}
\item all eigenvalues of $A$ are real and simple whenever
$b_2c_2,\dots,b_nc_n>0$ and $a_1,\dots,a_n\in\mathbb{R}$;
\item all eigenvalues of $A$ lie in the open right (resp., left)
half-plane whenever
 $b_2c_2,$ $\dots$, $b_nc_n<0$ and $a_1,\dots,a_n>0$ (resp., $<0$).
\end{itemize}
\end{problem}
\noindent{\bf Hint:}
Consider the rational function
$R(\lambda)\equiv\Delta_{n-1}(\lambda)/\Delta_n(\lambda)$
where $\Delta_n(\lambda)$ denotes the determinant of the matrix
 $\lambda I-A$ and  $\Delta_{n-1}(\lambda)$ denotes its principal
minor obtained by omitting its last row and column.
For this function, obtain a decomposition  of type
\eqref{eq55} or \eqref{eq56} and use the idea from
Remark~\ref{not75} and Problem~\ref{pro71}.
\begin{problem}\label{pro73}
Let $A:\mathfrak B\to\mathfrak B$ be a bounded linear operator
on a Banach space $\mathfrak B$, let $u\in\mathfrak B$ be a vector
in $\mathfrak B$, and let $\phi\in\mathfrak B^\prime$ be a bounded linear
functional acting on $\mathfrak B$ such that $\phi(u)\ne0$. We introduce
\begin{itemize}
\item the subspace
$\mathfrak B_1=\{x\in\mathfrak B:\phi(x)=0\}\subset\mathfrak B$,
\item the operator
$A_1:\mathfrak B_1\to\mathfrak B_1, \quad
A_1:x\mapsto Ax-\frac{\phi(Ax)}{\phi(u)}\,u$,
\item the vector $u_1\in\mathfrak B_1$ by $\quad u_1=Au-\frac{\phi(Au)}{\phi(u)}\,u$,
\item the functional
$\phi_1\in\mathfrak B_1^\prime,\quad \phi_1:x\mapsto\phi(Ax)$,
\end{itemize}
and consider the two analytic functions
\begin{equation*}
R(\lambda)\equiv\phi\left((\lambda I-A)^{-1}u\right),\quad
R_1(\lambda)\equiv\phi_1\left((\lambda I-A_1)^{-1}u_1\right).
\end{equation*}
Prove that
\begin{equation*}
-R(\lambda)=-\frac{s_0^2}{\lambda s_0-s_1-R_1(\lambda)},\quad\text{where}\quad
s_k\equiv\phi(A^ku)\quad (k=0,1).
\end{equation*}
Consequently, if $\dim\mathfrak B<\infty$, then $\IndC(R)=\sgn s_0+\IndC(R_1)$.
\end{problem}
\begin{problem}\label{pro74}
Prove that the projective transformations
\begin{equation}\label{eq57}
\Phi_\tau:\PR^1\to\PR^1,\quad
\Phi_\tau:\omega\mapsto\frac\omega{1-\tau\omega}\qquad(\tau\in\mathbb{R})
\end{equation}
form a {\it one-parameter group\/,} i.e.,
\begin{equation*}
\Phi_t\circ\Phi_\tau=\Phi_{t+\tau},\quad
\Phi_\tau^{-1}=\Phi_{-\tau}
\end{equation*}
\end{problem}
\begin{problem}\label{pro75} Using the notation of Problem~\ref{pro73},
consider
\begin{equation*}
A_t:\mathfrak B\to\mathfrak B,\quad
A_t:x\mapsto Ax+t\,\phi(x)\,u,\quad
R_t(\lambda)\equiv\phi\left((\lambda I-A_t)^{-1}u\right),
\end{equation*}
and let $\Phi_\tau$ be defined as in \eqref{eq57}. Prove that
$R_{t+\tau}=\Phi_{\tau}\circ R_t$.
Consequently, if $\dim\mathfrak B<\infty$, then
\begin{equation*}
\IndC(R_t)=\IndC(R_\tau)\quad(\forall t,\tau\in\mathbb{R}).
\end{equation*}
\end{problem}
\begin{problem}\label{pro76}
Theorem~\ref{thm71} dealt with the {\it left\/} action of the group
 $\SL,$ where the group acts on the {\it value\/} of the function $R$.
Prove that $\IndPR(R)$ is also invariant under the  {\it right\/}
action of the group $\SL$, where the group acts on the
{\it argument:\/}
\begin{equation*}
\IndPR(R\circ\Phi)=\IndPR(R),\quad\forall\Phi\in\SL .
\end{equation*}
\end{problem}
\begin{problem}\label{pro77}
If polynomials $f_0$ and $f_1$ are as in \S\,\ref{sec5}, 
then the function $R=f_1/f_0$ is odd, and Routh's algorithm, when it halts, yields
an expansion of $R$ into a continued fraction of the following form:
\begin{equation*}
R(\omega)=
\cfrac1{c_1\omega-
 \cfrac1{c_2\omega-
  \cfrac1{\ddots-
   \cfrac1{c_n\omega}}}} .
\end{equation*}
Prove that the following expansions are valid as well:
\begin{equation*}
R(\omega)=
\begin{cases}
\omega^{-1}\cdot
\cfrac1{c_1-
 \cfrac1{c_2\omega^2-
  \cfrac1{c_3-
   \cfrac1{\ddots-
    \cfrac1{c_{2k-1}}}}}}
\quad&(n=2k-1),\\
\omega\cdot
\cfrac1{c_1\omega^2-
 \cfrac1{c_2-
  \cfrac1{c_3\omega^2-
   \cfrac1{\ddots-
    \cfrac1{c_{2k}}}}}}
\quad&(n=2k).
\end{cases}
\end{equation*}
\end{problem}
\begin{problem}\label{pro78}
The expression
\begin{equation*}
\{R,z\}\equiv\frac{R^{\prime\prime\prime}(z)}{R^\prime(z)}-\frac32
\left[\frac{R^{\prime\prime}(z)}{R^\prime(z)}\right]^2
\end{equation*}
is called the {\it differential Schwarz invariant\/}, or
the {\it Schwarzian derivative\/} of the function $R$ (see \cite{19}).
Prove that
\begin{equation*}
\{\Phi\circ R,z\}=\{R,z\}\quad\forall\Phi\in\pmb{S}\pmb{L}(2,\mathbb{C}).
\end{equation*}
If $R(z)=s_0z+s_1z^2+s_2z^3+\cdots$, prove that
\begin{equation*}
\{R,0\}=\frac6{s_0^2}
\begin{vmatrix}
s_0 &s_1\\
s_1 &s_2
\end{vmatrix} .
\end{equation*}
\end{problem}


\section{Hermite-Biehler Theorem}\label{sec8}

Let a {\bf polynomial}
\begin{equation*}
p(z)=a_0z^n+a_1z^{n-1}+\dots+a_n\quad(a_0>0;\ a_1,\dots,a_n\in\mathbb{R})
\end{equation*}
be {\bf stable}. According to Theorem~3.1, the function $R\equiv f_1/f_0$,
with real polynomials $f_0$ and $f_1$ that are defined by
\begin{equation}\label{eq58}
i^{-n}p(i\omega)=f_0(\omega)-if_1(\omega)\qquad(\omega\in\mathbb{R}),
\end{equation}
is  {\it proper\/}: $\IndC(R)=\deg f_0=n$. Since the sum \eqref{eq17}
(or, equivalently, \eqref{eq47}) has no more than
$n$ terms $\pm1$, this equality is possible only if  the number of
terms is exactly $n$, and they are all equal to $+1$. Therefore,
the polynomial $f_0$ necessarily has $n$ distint real roots
 $\omega_1<\dots<\omega_n$, and since $\deg f_0=n$, it cannot
have additional (nonreal or multiple) roots. Hence $R$ splits
into {\it elementary fractions\/} as follows:
\begin{equation}\label{eq59}
R(z)=\sum\limits_{k=1}^n\frac{\alpha_k}{z-\omega_k},\quad
\alpha_k=\Res\limits_{\omega_k}(R) .
\end{equation}
By \eqref{eq16}, $\Ind_{\omega_k}(R)=\sgn\alpha_k$, and all indices
are equal to $+1$ in our case, so all residues $\alpha_k$ must be positive.
Thus
\begin{equation*}
\begin{matrix}
\dfrac{\im R(z)}{\im z}
&=-\sum\limits_{k=1}^n\dfrac{\alpha_k}{|z-\omega_k|^2}
<0\qquad
(\im z\neq0) , \\
\dfrac{d\,R(z)}{d\,z}
&=-\sum\limits_{k=1}^n\dfrac{\alpha_k}{(z-\omega_k)^2}
<0\qquad
(\im z=0) .
\end{matrix}
\end{equation*}

The function $-R$ therefore maps the upper half-plane
$\{z\colon\,\im z>0\}$ into itself, and monotonically increases between
its consecutive (real) poles. So, between any two consecutive roots
$\omega_{k-1}$, $\omega_k$ ($k=2,\dots,n$) of the denominator $f_0$
there must lie exactly one (simple) root of the numerator $f_1$.
Since $\deg f_1 \le n-1$, the polynomial $f_1$ cannot have any additional
roots.

Furthermore, the formula~\eqref{eq58} is {\it algebraic,\/} so
it does not matter that the argument $\omega$ was initially
assumed to be real. This formula may be re-written as
\begin{equation}\label{eq60}
i^{-n}p(z)=f_0(-iz)-if_1(-iz)\qquad(z\in\mathbb{C}) .
\end{equation}
Hence
\begin{equation*}
p(z)=0\,\implies\,
R(-iz)=-i\,\implies\,\frac{\im(-iz)}{\im(-i)}=\re z<0 .
\end{equation*}
{\bf The polynomial} $p$ {\bf is stable}! We now summarize this walk
{\rm\it<<there and back again>>}:
\begin{thm}\label{thm81}
Given a polynomial $p$,  let polynomials $f_0$ and $f_1$ be defined
by \eqref{eq58}, and let $R\equiv \left.f_1\right/f_0$. The following
conditions are equivalent:
\begin{enumerate}
\item the polynomial $p$ is stable;
\item the function $R$ is proper;
\item the function $R$ admits a representation of type \eqref{eq59},
with $\omega_k\in\mathbb{R}$ and $\alpha_k>0$ ($k=1,\dots,n$);
\item the function $-R$ maps the upper half-plane into itself;
\item the roots of the polynomials $f_0$, $f_1$ are real and simple,
 between any two consecutive roots of $f_0$ there is exactly one
root of $f_1$, and%
      \footnote{The last condition is a boring add-on. The pair
      $f_0$, $f_1$ should be normalized so that the ratio $f_1/f_0$
       be an decreasing rather than an increasing function over $\mathbb{R}$.}
\begin{equation}\label{eq61}
      \exists\,\omega\in\mathbb{R}\ :\
      f_1^\prime(\omega)\,f_0(\omega)-f_0^\prime(\omega)\,f_1(\omega)<0 .
\end{equation}
\end{enumerate}
\end{thm}
\begin{note}\label{not81}
The statement that conditions (1) and (5) are equivalent  is called
the {\it Hermite--Biehler\/} theorem. The two mathematicians obtained this result
simultaneously (1879) and independently - see \cite{1, 2, 3};
analogues for entire functions are given in  \cite{7}.
\end{note}
\begin{proof}
We already know that
\begin{equation*}
\begin{matrix}
(1)      &\Longleftrightarrow &(2)        &{}              &{} \\
\Uparrow &{}                  &\Downarrow &{}              &{} \\
(4)      &\Longleftarrow      &(3)        &\Longrightarrow &(5)
\end{matrix}
\end{equation*}
To prove the implication $(3)\Longleftarrow(5)$, note that the reality and
simplicity of the roots of $f_0$ imply the possibility of a decomposition of
type \eqref{eq59}. Next, if two consecutive residues $\alpha_{k-1}$ and
$\alpha_k$ were of different sign, then the interval  $(\omega_{k-1},\omega_k)$
would contain an {\it even\/} number of roots of $R$. All residues $\alpha_k$
are therefore of the same sign, namely positive, in view of \eqref{eq61}.
\end{proof}
\begin{note}\label{not82}
Analytic functions that map the upper half-plane into itself are well studied.
They play an important role in the spectral theory of self-adjoint operators
(see \cite{17, 18a, 18b}). The generic representation of such a function is
\begin{equation}\label{eq62}
F(z)=\alpha z+\beta+
\int\limits_{-\infty}^{+\infty}\frac{1+\omega z}{\omega-z}\,d\theta(\omega)
\qquad(\im z\ne0),
\end{equation}
where $\alpha \ge 0$, $\beta\in\mathbb{R}$, and $\theta(\omega)$ is a nondecreasing
function with finite limits  $\theta(\pm\infty)$. The function $\theta$ has only
a finite number of growth points  if and only if  $R=-F$ is a proper rational
function.
\end{note}
\begin{problem}\label{pro81}
Prove that the logarithmic derivative  of a polynomial $f$ with (not necessarily
real) roots $\omega_1,\dots,\omega_n$ satisfies \eqref{eq59},
where the residues  $\alpha_k$ are equal to the {\it multiplicities\/} of the roots
$\omega_k$.
\end{problem}
\begin{problem}\label{pro82}
Without recourse to \eqref{eq62}, prove that a real rational function $F$ maps
the upper half-plane into itself if and only if
\begin{equation*}
F(z)=\alpha z+\beta+\sum\limits_{k=1}^n\frac{\alpha_k}{\omega_k-z},\quad
\text{where}\ \alpha \ge 0,\ \beta\in\mathbb{R},\ \alpha_k>0,\ \omega_k\in\mathbb{R}.
\end{equation*}
\end{problem}
\begin{problem}\label{pro83}
Let $A$ be a bounded self-adjoint operator on a Hilbert space
 $\mathfrak H$, let $u\in\mathfrak H$, and consider
\begin{equation}\label{eq63}
R(\lambda)\equiv\left((\lambda I-A)^{-1}u,u\right) .
\end{equation}
Prove that
\begin{equation*}
\frac{\im R(\lambda)}{\im\lambda}=-\left\|(\lambda I-A)^{-1}u\right\|^2\quad
(\im\lambda\ne0)
\end{equation*}
and therefore  $-R$ maps the upper half-plane into itself.
\end{problem}
\begin{problem}\label{pro84}
Assuming that the operator of Problem \ref{pro83} is continuous,
prove that the corresponding function \eqref{eq63}  admits the decomposition
\begin{equation*}
R(\lambda)=\sum\limits_{k=1}^\infty\frac{|(u,e_k)|^2}{\lambda-\omega_k},
\end{equation*}
where $\{e_1,e_2,\dots\}$ is an orthonormal basis consisting of eigenvectors
of $A$, and  $\{\omega_1,\omega_2,\dots\}$ are the corresponding eigenvalues.
\end{problem}
\begin{problem}\label{pro85}
Let $R$ be a rational function that vanishes at infinity. Then it admits
a Laurent series representation that converges for sufficiently large $|z|$:
\begin{equation*}
R(z)=\frac{s_0}z+\frac{s_1}{z^2}+\frac{s_2}{z^3}+\dots
\end{equation*}
Prove that $R$ is proper if and only if the coefficients
$s_0,s_1,s_2,\dots$  satisfy
\begin{equation}\label{eq64}
s_k=\sum\limits_{j=1}^n\alpha_j\omega_j^k\qquad
\text{for some} \quad \alpha_1,\dots,\alpha_n>0, \; \; \;
\omega_1<\dots<\omega_n .
\end{equation}
\end{problem}
\begin{problem}\label{pro86}
If $R\in\mathcal R$, then
$-\dfrac{\im R(z)}{\im z}>\left|R^\prime(z)\right|$\quad $(\im z>0)$.
\end{problem}

Let $\gamma$ be a smooth curve lying in the upper half-plane
\begin{equation*}
\Pi\equiv\left\{\,z\in\mathbb{C}\,|\,\im z>0\,\right\}.
\end{equation*}
Define a {\it non-Euclidean\/} curve length $\gamma$ by the formula
$
s(\gamma)\equiv\int_\gamma (\im z)^{-1\,}|dz|.
$
This makes the half-plane $\Pi$ into the {\it Poincar\'e model\/} of the
Lobachevsky plane (see \cite{13, 23}). The geodesics of this plane
are half-circles with center on the real axis  $\mathbb{R}$, and
vertical rays. The Lobachevsky-Poincar\'e plane is shown on the front page
in a somewhat stylized form.

\begin{problem}\label{pro87}
If $R\in\mathcal R$, then the map
$-R\,:\,\Pi\to\Pi$ decreases curve lengths on $\Pi$:
\begin{equation*}
s(-R\circ\gamma)<s(\gamma) .
\end{equation*}
\end{problem}
\begin{problem}\label{pro88}
Let  $R\in\mathcal R$,
$\alpha,\beta,\gamma,\delta\in\mathbb{R}$. Then the equation
$R(z)=\dfrac{\alpha z+\beta}{\gamma z+\delta}$
\begin{itemize}
\item has only real solutions if
          $\alpha\delta-\beta\gamma \ge 0$;
\item has no more than one pair of complex conjugate solutions if
          $\alpha\delta-\beta\gamma<0$;
\end{itemize}
\end{problem}
\begin{problem}\label{pro89}
View the function $p(z)=e^z$ as being analogous to a polynomial.
This entire function has no roots in the {\it right\/} half-plane
(in fact, no roots whatsoever) so can be thought of as  <<stable>>.
Can one apply the results of Theorem~\ref{thm81}, at least to
some extent, to this function?  Without succumbing to premature
enthusiasm, consider also the function $p(z)=e^{-z}$.
\end{problem}


\section{Hankel forms}\label{sec9}

Given a real rational function%
\footnote{In contrast to \S\,\ref{sec7}, where a rational function was considered
as a map on the real projective line, here we consider it as an analytic function
of a complex variable.}
$R$, let us associate with it a {\it sesquilinear form\/}%
\footnote{$H\ \longleftarrow\ $
{\bf H}ermite, {\bf H}ankel, {\bf H}urwitz.}
$H$ defined on the (infinite-dimensional) complex linear space  $\mathcal P$ of
{\it all\/} polynomials by the formula
\begin{equation}\label{eq65}
H(x,y)\equiv\frac1{2\pi i}\,\oint\limits_\gamma R(\zeta)\,x(\zeta)\,
                                              \overline{y}(\zeta)\,d\zeta\qquad
                                               (x,y\in\mathcal P) .
\end{equation}
Here $\overline{y}(\zeta)\equiv\overline{y(\overline{\zeta})}$, and $\gamma$
is a positively oriented closed contour enclosing all poles of the function $R$.
By the Cauchy residue theorem,
\begin{equation}\label{eq66}
H(x,y)=\sum\Res\limits_{\omega_k}(Rx\overline{y})
\end{equation}
where the summation is over all poles of $R$.

A residue is especially easy to compute at a simple pole:
\begin{equation*}
f_0(\omega)=0,\ f_0^\prime(\omega)\ne0,\ f_1(\omega)\ne0\ \implies\
\Res\limits_\omega\left(\frac{f_1}{f_0}\right)=
\frac{f_1(\omega)}{f_0^\prime(\omega)} .
\end{equation*}
So, if all poles  $\omega_1,\dots,\omega_n$ of the function $R$
are real and simple, then
\begin{equation}\label{eq67}
H(x,x)=\sum\limits_{k=1}^n\alpha_k\left|x(\omega_k)\right|^2,\quad\text{where}
\quad\alpha_k\equiv\Res\limits_{\omega_k}(R)\in\mathbb{R} .
\end{equation}
We thus reduced the {\it Hermitian form\/} $H(x,x)$ to a sum of squares,
and formula \eqref{eq67} shows that the {\it rank\/} $\ran H$ of this
form (i.e., the total number of squares) is equal to $n$, the number of 
poles of the function $R$, whereas its {\it signature\/} $\sgn H$
(i.e., the difference between the number of positive and negative squares)
is equal to $\IndC(R)$ (Hermite 1856).

When the function $R$ has multiple or nonreal poles, one fails to
find such an expressive formula as \eqref{eq67}. However, the 
qualitative connection remains valid:
$\ran H = \deg f_0$,   $\sgn H = \IndC(R)$  (Hurwitz, 1895). 
Let us try to penetrate the essence of this phenomenon.

Given a polynomial $g\in\mathcal P$, consider the subspaces
\begin{equation*}
\mathcal P^g\equiv\{gu\,:\,u\in\mathcal P\},\qquad
\mathcal P_g\equiv\{v\in\mathcal P\,:\,\deg v<\deg g\} .
\end{equation*}
The {\it restrictions\/} of the form $H$ to these subspaces will be denoted
by $H^g$ and $H_g$, respectively. Long division of polynomials
 ($x=gu+v$, $\deg v<\deg g$) shows that%
\footnote{$\codim Y$ is the {\it codimension\/} of a subspace $Y$%
of a linear space $X$. It is equal to the dimension of a maximal
subspace of  $Z\subset X$ that intersects $Y$ trivially (according to
{\it Zorn's lemma}, such a subspace $Z$ always exists). If
$\dim X<\infty$, then $\codim Y=\dim X-\dim Y$.}
\begin{equation}\label{eq68}
\mathcal P=\mathcal P^g\oplus\mathcal P_g\quad\implies\quad
\codim\mathcal P^g=\dim\mathcal P_g=\deg g .
\end{equation}
\begin{note}\label{not91}
In addition to being a linear space, $\mathcal P$ is a
{\it commutative algebra\/}: its elements are multiplied according
to  well-known  rules  (see \cite{20, 21a, 21b, 21c}).
The subspace $\mathcal J=\mathcal P^g$ is an {\it ideal\/}
of the algebra $\mathcal P$:
\begin{equation*}
x\in\mathcal J, y\in\mathcal P\ \implies\ xy\in\mathcal J .
\end{equation*}
\end{note}

We now assign the <<inconvenient>> poles of the function $R=f_1/f_0$
to a polynomial $g$:
\begin{itemize}
\item    if $\omega$ is a nonreal pole of $R$ of order $\nu$ lying, say,
in the upper half plane, then it will be a root of $g$ of multiplicity
$\nu$ (then $\overline{\omega}$ will be a root of $\overline{g}$);
\item    if $\omega$ is a real pole of $R$ of order $\nu>1$, then it will
be a root of $g$ of  multiplicity
          $\lfloor\tfrac\nu2\rfloor\ge 1$.
\end{itemize}
As a result, the denominator $f_0$ splits (assuming that $f_0$ and $f_1$ have
no factors in common) into the product
\begin{equation}\label{eq69}
f_0=g\,\overline{g}\,h,
\end{equation}
where the roots of $h$ are all simple and coincide with real poles of $R$ of  {\it odd order\/}.

Now assume that $x,y\in\mathcal P^g$, $x=g\,x^g$, $y=g\,y^g$, and let
$R^g\equiv g\overline{g}R=f_1/h$. Then \eqref{eq65} turns into
\begin{equation}\label{eq70}
H^g(x,y)=
\frac1{2\pi i}\,\oint\limits_\gamma
R^g(\zeta)\,x^g(\zeta)\,\overline{y^g}(\zeta)\,d\zeta.
\end{equation}
All poles of the function $R^g$ are real and simple; therefore,
by the argument above,
\begin{equation}\label{eq71}
\ran H^g=\deg h,\qquad\sgn H^g=\IndC(R^g) .
\end{equation}
On the other hand, the polynomial $g\overline{g}$ is nonnegative
on $\mathbb{R}$, hence the functions $R^g$ and $R$ have equal signs
in the neighborhood of their common real poles, and hence
\begin{equation}\label{eq72}
\IndC(R^g)=\IndC(R) .
\end{equation}

We are now in a position to state and prove the main theorem of this section.
\begin{thm}\label{thm91}
The rank of the form $H$ is equal to the number of poles of the rational
function  $R$ (counted according to their multiplicities), and its signature
coincides with the Cauchy index $\IndC(R)$.
\end{thm}
\begin{proof}
The following two lemmata are necessary and sufficient for our proof.
\begin{lem}[rank]\label{lem91}
$
\ran H=\ran H^g+2\,\deg g=\deg f_0 .
$
\end{lem}
\begin{proof}
The second equality follows from~\eqref{eq69} and \eqref{eq71}.
We now prove the first.

By the theory of quadratic and Hermitian forms (see \cite{20, 21a, 21b, 21c}),
\begin{equation*}
\ran H=\codim\mathcal N,\quad\text{where}\quad
\mathcal N\equiv\{x\in\mathcal P\,:\,H(x,y)=0,\ \forall y\in\mathcal P\}
\end{equation*}
is a subspace called the {\it kernel\/} of $H$. To describe the kernel,
note that
\begin{equation*}
x\in\mathcal N\ \iff\
\oint\limits_\gamma R(\zeta)x(\zeta)\,\zeta^j\,d\zeta=0\quad
(j=0,1,2,\dots),
\end{equation*}
which means that $Rx$ is a polynomial, $x$ is divisible by the denominator
$f_0$, and $x\in\mathcal P^{f_0}$. Thus,
$\mathcal N=\mathcal P^{f_0}\ \implies\ \codim\mathcal N=\deg f_0$ .
\end{proof}

\begin{lem}[signature]\label{lem92}
$\sgn H=\sgn H^g$ .
\end{lem}
\begin{proof}
Denote the number of positive (resp., negative)  squares of the form $H$
by $\Pos H$  (resp., $\Neg H$). The number  $\Pos H$ coincides with the
dimension of a maximal subspace where the form  $H(x,x)$ is positive definite
(see \cite{20, 21a, 21b, 21c}); the same holds for  $\Pos H^g$, $\Neg H$, $\Neg H^g$.
Lemma \ref{lem91} implies that all these quantities are finite, and moreover,
\begin{equation}\label{eq73}
\Pos H+\Neg H=(\Pos H^g+\deg g)+(\Neg H^g+\deg g) .
\end{equation}

Let $\mathcal P^+$ be a subspace of $\mathcal P$ of dimension $\Pos H$, where
$H|_{\mathcal P^+}$ is positive definite. Then the restricted form $H^g$ is positive definite
on the subspace $\mathcal P^g\cap\mathcal P^+$, so that
\begin{equation}\label{eq74}
\dim(\mathcal P^g\cap\mathcal P^+)\le \Pos H^g .
\end{equation}
On the other hand,
\begin{equation}\label{eq75}
\dim\mathcal P^+\le \dim(\mathcal P^g\cap\mathcal P^+)+\codim\mathcal P^g .
\end{equation}
(Indeed, if $\mathcal Q$ is a maximal subspace of
$\mathcal P^+$ that intersects $\mathcal P^g$ trivially, then
$\dim\mathcal Q=\dim\mathcal P^+-\dim\left(\mathcal P^g\cap\mathcal P^+\right)$ and
$\dim\mathcal Q\le \codim\mathcal P^g$.) Since $\codim\mathcal P^g=\deg
g$, formul{\ae}  \eqref{eq74}--\eqref{eq75} imply
\begin{equation*}
\Pos H\le \Pos H^g+\deg g\quad\text{and,\/ analogously,}\quad
\Neg H\le \Neg H^g+\deg g,
\end{equation*}
but \eqref{eq73} implies that these must be equalities. Thus,
\begin{equation*}
\sgn H^g=\Pos H^g-\Neg H^g=\Pos H-\Neg H=\sgn H .
\end{equation*}
\end{proof}

This  concludes the proof of Theorem~\ref{thm91} as well: the statement about
the rank follows from Lemma~\ref{lem91}, and the statement about the signature
follows from Lemma~\ref{lem92} and from equalities \eqref{eq71}--\eqref{eq72} .
\end{proof}

\begin{note}\label{not92}
We now understand the influence of nonreal and multiple roots
that were collected into $g$: they give the form $H$ some <<ballast>>
consisting of an equal number ($\deg g$) of positive and negative
squares, which brings up the rank of $H$.  These rather nondescript
squares get filtered out when $H$ is restricted to the ideal
$\mathcal P^g$ of the algebra  $\mathcal P$.
\end{note}

\begin{note}\label{not93}
Theorem \ref{thm91} deals with $\IndC(R)$ -- not with $\IndPR(R)$, as in
 \S\,\ref{sec7}. However, if $\deg f_1\le\deg f_0$, the two indices coincide.
\end{note}

So far, it was convenient%
\footnote{Does the reader see why?}
to consider the form $H$ on a {\it complex infinite-dimensional\/}
space. We take a more practical position now.

First of all, the form  $H$ is real:
$\overline{H(x,y)}=H(\overline{x},\overline{y})$. This implies that, for
real polynomials $u$ and $v$,  the equality $H(u+iv,u+iv)=H(u,u)+H(v,v)$ holds.
So, without loss of generality, $\mathcal P$ can be assumed to be a {\it real\/}
space (or algebra).

Secondly, denote for simplicity $f\equiv f_0$ and recall that
 $\mathcal P=\mathcal P_f\oplus\mathcal P^f$ (relation \eqref{eq68})
as well as $\mathcal P^f=\mathcal N$ (the proof of Lemma~\ref{lem91}).
This implies that
\begin{equation}\label{eq76}
H(x_f+x^f,y_f+y^f)=H_f(x_f,y_f)\qquad
(x_f,y_f\in\mathcal P_f;\ x^f,y^f\in\mathcal P^f).
\end{equation}
\begin{lem}\label{lem93}
$\ran H=\ran H_f$, $\sgn H=\sgn H_f$
\end{lem}
\begin{proof}
With $\mathcal P^+$ as in the proof of Lemma~\ref{lem92},
formula $\eqref{eq76}$ yields
\begin{equation*}
x\in\mathcal P^+,\ x\ne0\ \implies\ H(x,x)=H_f(x_f,x_f)>0\ \implies\
x_f\ne0.
\end{equation*}
So, the projector  $x\mapsto x_f$ to the first component of the direct sum%
\footnote{Simply, this is division of the polynomial $x$ with quotient $f$
and remainder $x_f$.}
$\mathcal P_f\oplus\mathcal P^f$ is injective on $\mathcal P^+$, hence does
not decrease its dimension. Therefore $\Pos H_f\ge\Pos H$, while the reverse
inequality is obvious. Analogously, $\Neg H_f=\Neg H$.
\end{proof}

Thus, we can consider the form $H$ on the  {\it finite-dimensional\/}
space $\mathcal P_f$. However, $\mathcal P_f$, unlike $\mathcal P$, is not
an algebra. On the other hand, the form $H_f$ is necessarily {\it nondegenerate}:
\begin{equation*}
x\in\mathcal P_f,\ H_f(x,y)=0,\ \forall y\in\mathcal P_f\quad
\implies\quad x=0 .
\end{equation*}

Finally, consider a basis in $\mathcal P_f$ and write the form
$H_f$ in its canonical form.

Let the basis consist of monomials $\zeta^j$ ($j=0,1,\dots,n-1$;
$n\equiv\deg f$, $f\equiv f_0$).
Then \eqref{eq65} shows that
\begin{equation}\label{eq77}
x(\zeta)=\sum\limits_{j=0}^{n-1}\xi_j\zeta^j\quad\implies\quad
H_f(x,x)=\sum\limits_{i,j=0}^{n-1}s_{i+j}\xi_i\xi_j,
\end{equation}
where
\begin{equation}\label{eq78}
s_k\equiv
\frac1{2\pi i}\oint\limits_\gamma R(\zeta)\,\zeta^k\,d\zeta\quad
(k=0,1,2,\dots) .
\end{equation}
If the rational function $R$ is expanded into a series
\begin{equation}\label{eq79}
R(\zeta)=
s_{-m}\zeta^{m-1}+\dots+\frac{s_0}\zeta+\frac{s_1}{\zeta^2}+\cdots\qquad
(m=\deg f_1-\deg f_0)
\end{equation}
(which converges absolutely for large values of $|\zeta|$ if $m\le 1$),
then, substituting \eqref{eq79} into \eqref{eq78} and integrating
term by term, we verify that the coefficients $s_k$ ($k\ge 0$) in \eqref{eq79}
satisfy conditions~\eqref{eq78}.

An unusual property of \eqref{eq77} is that each of its coefficients
$s_{i+j}$ depends only on the {\it sum\/} of its indices.
Quadratic forms of this type are called {\it Hankel forms.}

The next theorem  adds to our already large collection of statements
that characterize rational functions from the class $\mathcal R$:
\begin{thm}\label{thm92} Let
\begin{equation*}
R(\zeta)=\frac{f_1(\zeta)}{f_0(\zeta)}=
\frac{s_0}\zeta+\frac{s_1}{\zeta^2}+\cdots\quad
(\deg f_1<\deg f_0=n).
\end{equation*}
Then $R\in\mathcal R$ if and only if
\begin{equation}\label{eq80}
s_0>0,\quad
\vmatrix
s_0 &s_1\\
s_1 &s_2
\endvmatrix>0,\quad
\vmatrix
s_0     &\dots &s_{n-1} \\
\vdots   &\ddots &\vdots   \\
s_{n-1} &\dots &s_{2n-2}
\endvmatrix>0.
\end{equation}
\end{thm}
\begin{proof}
The condition $R\in\mathcal R$ is equivalent to the positive definiteness
of the form~\eqref{eq77}, which is by Sylvester's criterion equivalent
to \eqref{eq80}.
\end{proof}
\begin{note}\label{not94}
It would be interesting to find out how the determinants in \eqref{eq80}
can be expressed directly in terms
of the coefficients of the polynomials $f_1$ and $f_0$. We will take up
this question in the next section.
\end{note}

\begin{problem}\label{pro91}
Prove that  $\mathcal P$ is a  {\it principal ideal domain.\/}
The latter means  \cite{20, 21a, 21b, 21c} that any ideal
$\mathcal J\subset\mathcal P$ (see Remark \ref{not91})
is generated by some polynomial  $g\in\mathcal P$:
$\mathcal J=\mathcal P^g$.
\end{problem}

\begin{problem}\label{pro92}
Let $f,g\in\mathcal P$ and let $d=\CMD(f,g)$, $k=\CMM(f,g)$.
Prove that $\mathcal P^k=\mathcal P^f\cap\mathcal P^g$ and
$\mathcal P^d=\mathcal P^f+\mathcal P^g$.
This explains the real meaning and role of the terms
<<greatest common factor>> and <<least common multiple>>.
\end{problem}

\begin{problem}\label{pro93}
Under the assumptions of Problem \ref{thm92}, suppose  the following
inequalities hold as well:
\begin{equation}\label{eq81}
s_1>0,\quad
\begin{vmatrix}
s_1 &s_2\\
s_2 &s_3
\end{vmatrix}
>0,\quad
\begin{vmatrix}
s_1     &\dots &s_n     \\
\vdots   &\ddots &\vdots   \\
s_n     &\dots &s_{2n-1}
\end{vmatrix}
>0 .
\end{equation}
What additional properties of the function $R$ follow?
\end{problem}

\begin{problem}\label{pro94}
Prove the {\it Borchardt-Jacobi Theorem:\/} Let $f$ be a real polynomial
of degree $n$ with (complex) roots $\lambda_1$, \dots,
$\lambda_n$, and let%
\footnote{The quantities being defined here are called Newton sums.
They are {\it symmetric\/} functions of the roots  $\lambda_j$;
therefore they can be found without knowing the actual roots
(see \cite{7}).}
\begin{equation*}
s_k=\lambda_1^k+\dots+\lambda_n^k\quad(k=0,1,2,\dots) .
\end{equation*}
Then the number of positive squares of the form
$\sum\nolimits_{i,j=0}^{n-1}s_{i+j}\xi_i\xi_j$ coincides with the
number of {\em distinct\/} roots of the polynomial $f$, and
the number of its negative squares with the number of distinct
complex conjugate pairs of roots.
\end{problem}
\noindent{\bf Hint:} One can apply the results of this and previous
Sections to
the logarithmic derivative of the polynomial  $f$, but a direct
argument is also possible.

\begin{problem}\label{pro95}
Derive {\it Newton's identities\/} that connect the
Newton sums  $s_0$, $s_1$, $\dots$ with the coefficients
 $a_0$, $a_1$, $\dots$, $a_n$ of the polynomial $f$ from
the previous problem.
\end{problem}
\noindent{\bf Hint:} Prove the relation
\begin{equation*}
\frac{na_0\zeta^{n-1}+\dots+a_{n-1}}{a_0\zeta^n+\dots+a_n}=
\frac{s_0}{\zeta}+\frac{s_1}{\zeta^2}+\frac{s_2}{\zeta^3}+\cdots
\end{equation*}
and make use of it.

\begin{problem}\label{pro96}
As in Problem~\ref{pro73}, let
$R(\lambda)\equiv\phi\left((\lambda I-A)^{-1}u\right)$,
where $A\ :\ \mathfrak B\to\mathfrak B$ is a bounded linear operator
on a Banach space $\mathfrak B$, $u\in\mathfrak B$,
$\phi\in\mathfrak B^\prime$. Prove that
\begin{equation*}
R(\lambda)=
\frac{s_0}{\lambda}+\frac{s_1}{\lambda^2}+\frac{s_2}{\lambda^3}+\cdots,\quad
\text{where}\quad s_k=\phi\left(A^ku\right)\quad(k=0,1,2,\dots).
\end{equation*}
\end{problem}

\begin{problem}\label{pro97}
As in problem \ref{pro83}, let
$R(\lambda)\equiv\left((\lambda I-A)^{-1}u,u\right)$,
where $A\ :\mathfrak H\to\mathfrak H$ is a bounded self-adjoint operator
on a Hilbert space $\mathfrak H$. Prove that the Hankel forms
\begin{equation*}
\sum\limits_{i,j=0}^{n-1}s_{i+j}\xi_i\xi_j\quad(n=1,2,\dots),\quad
\text{where}\quad s_k=\left(A^ku,u\right)\quad(k=0,1,\dots)
\end{equation*}
are nonnegative definite. In which case are they positive definite?
\end{problem}

\begin{problem}\label{pro98}
Let $A$ be a real  $n\times n$-matrix whose elements and minors are
all positive%
\footnote{Such matrices are called totally positive;
for details, see \cite{22}.}. Let
$s_k$ ($k=0,1,2,\dots$) be the  $(1,1)$ entry of the $k$th power $A^k$
of the matrix $A$.  Prove that inequalities \eqref{eq80} and
\eqref{eq81} hold.
\end{problem}
\noindent{\bf Hint:} Find formul{\ae} connecting the determinants in
\eqref{eq80}-\eqref{eq81} with the minors of  $A$.


\section{Li\'enard-Chipart criterion}\label{sec10}

We now attempt to combine the ideas of the preceding three sections.

The rank and the signature of a Hankel form are {\it not\/} projective invariants,
since they do not account for a possible pole of the function
$R=f_1/f_0$ at the point $\infty$. However, if $\deg f_1\leq\deg f_0$,
then $R(\infty)\ne\infty$ and $\IndPR(R)=\IndC(R)$. Moreover, if, as in
\S\,\ref{sec7},
\begin{equation*}
\Phi\in\SL,\quad
\Phi\,:\,\omega\mapsto\frac{\alpha\omega+\beta}{\gamma\omega+\delta}\quad
(\alpha\delta-\beta\gamma=1),
\end{equation*}
then
\begin{equation*}
(\Phi\circ R)(z)=
\frac{\alpha f_1(z)+\beta f_0(z)}{\gamma f_1(z)+\delta f_0(z)},
\end{equation*}
and
\begin{equation*}
\deg(\alpha f_1+\beta f_0)=
\deg(\gamma f_1+\delta f_0)=
\max\left\{\deg f_1,\deg f_0\right\}
\end{equation*}
for all $\Phi\in\SL$, except for exactly two elements: $\Phi_0$ when $\Phi_0\circ R$ has a zero at the
point  $z=\infty$, and  $\Phi_\infty$ if $z=\infty$ is a pole of $\Phi_\infty\circ R$.
\begin{lem}\label{lem10.1}
If
\begin{equation}\label{eq82}
R(z)\equiv\frac{b_0z^n+b_1z^{n-1}+\dots+b_n}{c_0z^n+c_1z^{n-1}+\dots+c_n}=
s_{-1}+\frac{s_0}z+\frac{s_1}{z^2}+\cdots\quad(c_0\ne0),
\end{equation}
then%
\footnote{Lo and behold the Hurwitz matrix!}
\begin{equation}\label{eq83}
\begin{array}{ll}
&\nabla_{2k}\equiv
\begin{vmatrix}
c_0    &c_1    &c_2    &\hdots &c_{k-1}  &c_k     &\hdots &c_{2k-1} \\
b_0    &b_1    &b_2    &\hdots &b_{k-1}  &b_k     &\hdots &b_{2k-1} \\
0      &c_0    &c_1    &\hdots &c_{k-2}  &c_{k-1} &\hdots &c_{2k-2} \\
0      &b_0    &b_1    &\hdots &b_{k-2}  &b_{k-1} &\hdots &b_{2k-2} \\
\vdots & \vdots & \vdots & \ddots & \vdots & \vdots &  \ddots & \vdots  \\
0      &0      &0      &\hdots &c_0      &c_1     &\hdots &c_k      \\
0      &0      &0      &\hdots &b_0      &b_1     &\hdots &b_k
\end{vmatrix}=
c_0^{2k}\det\left[s_{i+j}\right]_0^{k-1}\\
&(k=1,2,\dots,n;\quad\text{as in \S\,\ref{sec6}, we set}\quad
c_j=b_j=0\ \text{for}\ j>n).
\end{array}
\end{equation}
\end{lem}
\begin{proof}
First interchange the rows of the determinant in \eqref{eq83} to obtain
\begin{equation}\label{eq84}
\nabla_{2k}=
\begin{vmatrix}
c_0    &c_1    &c_2    &\hdots &c_{k-1} &\tvr  &c_k     &\hdots &c_{2k-1}\\
0      &c_0    &c_1    &\hdots &c_{k-2} &\tvr  &c_{k-1} &\hdots &c_{2k-2}\\
\vdots & \vdots & \vdots & \ddots & \vdots & & \vdots &  \ddots & \vdots \\
0      &0      &0      &\hdots &c_0     &\surd &c_1     &\hdots &c_k     \\
0      &0      &0      &\hdots &b_0     &\surd &b_1     &\hdots &b_k     \\
\vdots & \vdots & \vdots & \ddots & \vdots & & \vdots &  \ddots & \vdots \\
0      &b_0    &b_1    &\hdots &b_{k-2} &\tvr  &b_{k-1} &\hdots &b_{2k-2}\\
b_0    &b_1    &b_2    &\hdots &b_{k-1} &\tvr  &b_k     &\hdots &b_{2k-1}
\end{vmatrix}.
\end{equation}
This does not change the sign of $\nabla_{2k}$. Indeed, lower the  $k$th and $(k{+}1)$st rows%
\footnote{They are marked by $\surd$ in \eqref{eq84}.} to their initial
positions in \eqref{eq83}. This will require an  {\it even\/} number of
transpositions. The next pair of rows will then meet, the lowering
operation will be applied to them, and so on.

Now let us establish a connection between $b,c$ and $s$.
Multiplying \eqref{eq82} by the denominator and equating coefficients,
we get:
\begin{equation}\label{eq85}
b_j=\sum_{i=0}^j c_{j-i}s_{i-1}\qquad(j=0,1,2,\dots).
\end{equation}
These formul{\ae} suggest by themselves what to do next, namely, to
eliminate the entries in the lower left corner of the determinant \eqref{eq84}.
From each  $k+j$th row
($j=1,2,\dots,k$) subtract rows $k-j+1$, $k-j+2$, \dots, $k$ multiplied by
$s_{-1}$, $s_0$, $\dots$, $s_{j-2}$, respectively. As a result, we get zeros
down and to the left, and the entry
\begin{multline*}
d_{ij}\equiv b_{i+j-1}-c_{i+j-1}s_{-1}-\dots-c_is_{j-2}=
c_{i-1}s_{j-1}+\dots+c_0s_{i+j-2}\\
(i,j=1,2,\dots,k)
\end{multline*}
in position $(k+i,k+j)$. In matrix form, this will look as follows:
\begin{equation*}
\begin{vmatrix}
d_{11} &d_{12} &\dots &d_{1k}\\
d_{21} &d_{22} &\dots &d_{2k}\\
\vdots & \vdots & \ddots & \vdots \\
d_{k1} &d_{k2} &\dots &d_{kk}
\end{vmatrix}=
\begin{vmatrix}
s_0     &s_1   &\dots &s_{k-1} \\
s_1     &s_2   &\dots &s_k     \\
\vdots & \vdots & \ddots & \vdots \\
s_{k-1} &s_k   &\dots &s_{2k-2}
\end{vmatrix}\,\cdot\,
\begin{vmatrix}
c_0   &c_1   &\dots &c_{k-1}\\
0     &c_0   &\dots &c_{k-2}\\
\vdots & \vdots & \ddots & \vdots \\
0     &0     &\dots &c_0
\end{vmatrix} .
\end{equation*}
\end{proof}

\begin{note}\label{not10.1}
The quantities $\nabla_{2k}=\nabla_{2k}(R)$ are invariants under the action
of the group $\SL$, since the structure of determinants \eqref{eq83} implies
that
$\nabla_{2k}(R+d)=\nabla_{2k}(R)$ ($d=\text{const}$) and
$\nabla_{2k}(-\frac1R)=\nabla_{2k}(R)$ due to Lemmata 7.1--7.3.
\end{note}

\begin{note}\label{not10.2} $\nabla_{2n}$ coincides, up to a sign,
with the {\it resultant\/} of the polynomials $f_1$ and $f_0$. Hence
\begin{equation*}
\nabla_{2n}\ne0\quad\iff\quad\CMD(f_0,f_1)=1
\end{equation*}
More about the resultant and other  {\it symmetric polynomials\/} can be
found in \cite{7}; this topic is worth studying {\it per se,\/} but
we do not need to invoke external results to justify the fact we just stated.
Indeed, the inequality  $\nabla_{2n}\ne0$ together with Lemma \ref{lem10.1}
imply that the rank of the corresponding Hankel form $H$, which coincides
with the total number of the poles of the function $f_1/f_0$, is at least
 $\deg f_0$, hence the fraction $f_1/f_0$ is in lowest terms.
\end{note}

\begin{lem}\label{lem10.2}
With the notation of Lemma \ref{lem10.1},
\begin{equation*}
\nabla_2>0,\ \nabla_4>0,\ \dots,\ \nabla_{2n}>0\quad\implies\quad
\frac{\im R(z)}{\im z}<0\quad(\im z\ne0).
\end{equation*}
If $n$ is the number of poles of the rational function $R$, then the converse
holds as well.
\end{lem}
\begin{proof}
If $b_0=0$, apply Lemma \ref{lem10.1}, Theorem \ref{thm92}, and
Theorem \ref{thm81}. If $b_0\ne0$, apply the same results to the
function $R(z)-s_{-1}$.
\end{proof}

We already noticed that  the determinants $\nabla_{2k}$ ($k=1,2,\dots,n$) are
the even leading principal minors  of the $2n\times2n$  matrix
built exactly as the Hurwitz matrix $\mathcal H_p$ from \S\,6.
This is not just a superficial similarity. Let us go back  to the
main object of our study, viz., the real polynomial
\begin{equation}\label{eq86}
p(z)\equiv a_0z^n+a_1z^{n-1}+\dots+a_n\qquad(a_0>0)
\end{equation}
and let us rewrite it as
\begin{equation*}
p(z)=g_0(z^2)+zg_1(z^2).
\end{equation*}
The polynomials $g_0$ and $g_1$ so defined are very much reminiscent
of   $f_0$ and $f_1$, which first appeared in \S\,\ref{sec3} and were studied in
detail in \S\S\,\ref{sec5}, \ref{sec6}. If the degree $n$ is odd $n=2m+1$, then
\begin{equation*}
\begin{array}{ll}
&g_0(w)=a_1w^m+a_3w^{m-1}+\dots+a_{2m+1},\\
&g_1(w)=a_0w^m+a_2w^{m-1}+\dots+a_{2m};
\end{array}
\end{equation*}
if  $n=2m$, then
\begin{equation*}
\begin{array}{ll}
&g_0(w)=a_0w^m+a_2w^{m-1}+\dots+a_{2m},\\
&g_1(w)=a_1w^{m-1}+a_3w^{m-2}+\dots+a_{2m-1}.
\end{array}
\end{equation*}

In either case,  $\deg g_0\ge \deg g_1$, and the minors
$\nabla_{2k}=\nabla_{2k}(g_1/g_0)$ and the minors $\eta_k$
of the Hurwitz matrix $\mathcal H_p$ (see \S\,\ref{sec6}) are connected  as follows%
\footnote{Note that in both cases the {\it last\/} $\nabla_{2m}$
corresponds to the {\it second-last\/} $\eta_{n-1}$, which, as we know
from \S\,\ref{sec5},  <<guards>> the imaginary axis!}:
\begin{equation}\label{eq87}
\nabla_{2k}=
\begin{cases}
   \eta_{2k}  \quad(k=1,\dots,m),\qquad&\text{if}\quad n=2m+1\\
a_0\eta_{2k-1}\quad(k=1,\dots,m),\qquad&\text{if}\quad n=2m.
\end{cases}
\end{equation}

If, as usual, $R\equiv g_1/g_0$, then
\begin{equation}\label{eq88}
p(z)=0\quad\iff\quad R(z^2)=-\frac1z .
\end{equation}

We are now ready to prove the Li\'enard-Chipart criterion, which was
mentioned in \S\,\ref{sec6}:

\begin{thm}[Li\'enard--Chipart]\label{thm10.1}
A polynomial \eqref{eq86} is stable if and only if
\begin{align}
&a_n>0,\ a_{n-1}>0,\ a_{n-2}>0,\ \dots\label{eq89}\\
&\eta_{n-1}>0,\ \eta_{n-3}>0,\ \eta_{n-5}>0,\ \dots\label{eq90}
\end{align}
\end{thm}

\begin{proof}
Necessity: follows from the Stodola condition and the Hurwitz criterion.

\noindent
Sufficiency:
1. Since the minors are positive, Lemma \ref{lem10.2} and formula \eqref{eq87}
imply
\begin{equation*}
\frac{\im R(z)}{\im z}<0\qquad(\im z\ne0)
\end{equation*}

So, if $p(z)=0$, $\im z\ne0$, then \eqref{eq88} and \eqref{eq90} give
\begin{equation*}
\frac{\im R(z^2)}{\im z^2}=\frac{\im z}{|z|^2\im z^2}=\frac1{2|z|^2\re z}<0.
\end{equation*}
Thus, by  \eqref{eq90}, the roots of the polynomial $p$ lie in the union of the open
left-hand plane and the real positive half-line.

\noindent
2. Since the coefficients are positive, the polynomial $p$ cannot have
any roots on the nonnegative side of the real axis.
\end{proof}

\begin{note}\label{not10.3} If we drop \eqref{eq89} but keep \eqref{eq90} in
Theorem \ref{thm10.1}, then some roots of $p$ may cross over into the right
half-plane, but in that case {\it they must stay on the positive real half-line.\/}
From the stability point of view, this behavior is very interesting,
since it means that the fixed point  does not bifurcate into a limit cycle (these bifurcations are described in detail in \cite{24}).
\end{note}

\begin{problem}\label{pro10.1}
For a polynomial discussed in  Remark~\ref{not10.3} whose  roots with positive
real part cannot leave the real axis, it is natural to expect that these roots
stay simple (a simple root cannot leave the real axis other than by coalescing
with another root). Prove that this is indeed so, i.e., that the positive
roots of a polynomial satisfying  \eqref{eq90} are simple.
\end{problem}


\section*{Afterword}

Half a page is still left -- it would be sinful to leave it
blank. Let us draw conclusions.

In stability theory, the Routh-Hurwitz problem, which was
considered in these <<lecture notes>> from many points of view,
certainly does not play a role commensurate with the
attention we devoted to it; more precisely, it does not yet
play a role it is destined for. Destined by whom?
I don't know. Still, I believe that intrinsically beautiful
mathematical constructs must be necessarily connected to
the understanding and explanation of the real world that
surrounds us. If you wish, you may label this a mathematical
religion of sorts; I think many mathematicians,  perhaps most,
are such conscious or unconscious believers. It is impossible
to imagine that the unreal world studied by mathematics has
been created by human intellect; mathematicians do not invent
theorems and theories but {\it discover\/} them.  And this ideal
world of mathematical constructions always turns out to parallel
the <<real>> one; remarkably, its researchers are driven neither
by  <<logic>>, as laymen think, nor by <<applied>> needs, nor
even by acquired experience and knowledge, both certainly
indispensable. They are driven by an irrational, strange  intuition
that lets them feel that intrinsic beauty and harmony, just as
our senses can feel warmth and determine its source.
It is true that the mathematician Hurwitz was <<handed>> a problem by
the turbine engineer Stodola, but Hurwitz took up and solved that
problem not to help Stodola build his turbines. Well, he would
have not taken it up just for that.  Such is indeed the relationship
between mathematics and its various <<applications>>:  the latter are
sources of problems for the former, and the origin of these problems
is a certain {\it a priori\/} guarantee of the harmony and beauty to be
found there, and of the progress in mathematics that their solutions
must bring about.

Returning to the beginning of our scholia, let me venture an opinion:
the sad fact that the <<Routh-Hurwitz problem>>  is not much in demand
does not mean that it is only of academic interest. This should rather mean
that there is a hidden door behind which there may be lots of interesting
stuff. Just as in the instructive story  of ``The golden key''%
\footnote{"The golden key" by Alexei Tolstoi, a famous Russian adaptation of
the book "The adventures of Pinocchio" by Carlo Collodi [translators' remark].}. :--)

\def\cprime{$'$} \def\cprime{$'$}

\smallskip

\noindent{\it Comments.}%
\footnote{Many of the author's original references, in Russian or translated into Russian, were replaced by the corresponding references in English [translators' remark].}
The books \cite{1,2, 3} offer  a thorough treatment of the subject.
The book \cite{2} is most elementary and detailed. The monograph \cite{1}
remains one of the best matrix theory books in the world.
The paper \cite{3} contains an exhaustive review of classsical works
on stable polynomials. The textbook \cite{4} states the Hurwitz theorem
and the <<Mikhailov criterium>>. However, the relevant sections of this overall
good book are in my opinion not quite satisfactorily, and it is better to use~\cite{2}.

The monographs~\cite{6, 7} are devoted to the Routh-Hurwitz problem
for entire functions.

The books~\cite{9, 8, 10, 11a, 11b} of G. Polya mentioned in the Introduction are not
directly related  to our topic of stable polynomials, but their reading is useful
for every beginning mathematician.

The book~\cite{12} is outdated and is written <<for dummies>>. However,
the now forgotten language {\sc APL}, which was created by Kenneth Iverson not
exactly as a practical programming language but rather  as a notation system
for mathematical algorithms,  is {\it per se\/} interesting to a mathematician.


\begin{thebibliography}{10}

\bibitem{6}
N.~I. Ahiezer and M.~Krein.
\newblock {\em Some questions in the theory of moments}.
\newblock translated by W. Fleming and D. Prill. Translations of Mathematical
  Monographs, Vol. 2. American Mathematical Society, Providence, R.I., 1962.

\bibitem{17}
N.~I. Akhiezer and I.~M. Glazman.
\newblock {\em Theory of linear operators in {H}ilbert space. {V}ol. {I}},
  volume~9 of {\em Monographs and Studies in Mathematics}.
\newblock Pitman (Advanced Publishing Program), Boston, Mass., 1981.
\newblock Translated from the third Russian edition by E. R. Dawson,
  Translation edited by W. N. Everitt.

\bibitem{15}
V.~I. Arnol{\cprime}d.
\newblock {\em Dopolnitelnye glavy teorii obyknovennykh differentsialnykh
  uravnenii}.
\newblock ``Nauka'', Moscow, 1978.


 

\bibitem{21a}
N.~Bourbaki.
\newblock {\em \'{E}l\'ements de math\'ematique. {XIV}. {P}remi\`ere partie:
  {L}es structures fondamentales de l'analyse. {L}ivre {II}: {A}lg\`ebre.
  {C}hapitre {VI}: {G}roupes et corps ordonn\'es. {C}hapitre {VII}: {M}odules
  sur les anneaux principaux}.
\newblock Actualit\'es Sci. Ind., no. 1179. Hermann et Cie, Paris, 1952.

\bibitem{21b}
N.~Bourbaki.
\newblock {\em \'{E}l\'ements de math\'ematique. 23. {P}remi\`ere partie: {L}es
  structures fondamentales de l'analyse. {L}ivre {II}: {A}lg\`ebre. {C}hapitre
  8: {M}odules et anneaux semi-simples}.
\newblock Actualit\'es Sci. Ind. no. 1261. Hermann, Paris, 1958.

\bibitem{21c}
N.~Bourbaki.
\newblock {\em \'{E}l\'ements de math\'ematique. {P}remi\`ere partie: {L}es
  structures fondamentales de l'analyse. {L}ivre {II}: {A}lg\`ebre. {C}hapitre
  9: {F}ormes sesquilin\'eaires et formes quadratiques}.
\newblock Actualit\'es Sci. Ind. no. 1272. Hermann, Paris, 1959.

\bibitem{7}
N.~G. {\v{C}}ebotarev and N.~N. Me{\u\i}man.
\newblock The {R}outh-{H}urwitz problem for polynomials and entire functions.
  {R}eal quasipolynomials with {$r=3$}, {$s=1$}.  (Russian)
\newblock {\em Trudy Mat. Inst. Steklov.}, 26:331, 1949.
\newblock Appendix by G. S. Barhin and A. N. Hovanski\u\i.

\bibitem{4}
B.~P. Demidovi{\v{c}}.
\newblock {\em Lektsii po matematicheskoi teorii ustoichivosti}.
\newblock Izdat. ``Nauka'', Moscow, 1967.

\bibitem{13}
B.~A. Dubrovin, A.~T. Fomenko, and S.~P. Novikov.
\newblock {\em Modern geometry---methods and applications. {P}art {I}},
  volume~93 of {\em Graduate Texts in Mathematics}.
\newblock Springer-Verlag, New York, second edition, 1992.
\newblock The geometry of surfaces, transformation groups, and fields,
  Translated from the Russian by Robert G. Burns.

\bibitem{22}
F.~P. Gantmacher and M.~G. Krein.
\newblock {\em Oscillation matrices and kernels and small vibrations of
  mechanical systems}.
\newblock AMS Chelsea Publishing, Providence, RI, revised edition, 2002.
\newblock Translation based on the 1941 Russian original, Edited and with a
  preface by Alex Eremenko.

\bibitem{1}
F.~R. Gantmacher.
\newblock {\em The theory of matrices. {V}ols. 1, 2}.
\newblock Translated by K. A. Hirsch. Chelsea Publishing Co., New York, 1959.



\bibitem{12}
L.~Gilman and A.~J. Rose.
\newblock {\em APL -- an interactive approach}.
\newblock Third Edition. John Wiley and Sons. 1984.


\bibitem{14}
G.~B. Gurevi{\v{c}}.
\newblock {\em Osnovy {T}eorii {A}lgebrai\v ceskih {I}nvariantov}.
\newblock OGIZ, Moscow-Leningrad, 1948.

\bibitem{23}
D.~Hilbert and S.~Cohn-Vossen.
\newblock {\em Geometry and the imagination}.
\newblock Chelsea Publishing Company, New York, N. Y., 1952.
\newblock Translated by P. Nem\'enyi.

\bibitem{19}
A.~Hurwitz and R.~Courant.
\newblock {\em Vorlesungen \"uber allgemeine {F}unktionentheorie und
  elliptische {F}unktionen}.
\newblock Interscience Publishers, Inc., New York, 1944.

\bibitem{18a}
I.~S. Kac and M.~G. Kre\u\i n. {\em $R$-functions  -- analytic functions mapping the
upper halfplane into itself.}
\newblock Translations Am. Math. Soc., ser.~2, 103 : 1--18, 1974.

 \bibitem{18b}
I.~S. Kac and M.~G. Kre\u\i n. {\em On the spectral functions of the string.}
\newblock Translations Am. Math. Soc., ser.~2, 103 : 19--102, 1974.


\bibitem{3}
M.~G. Kre{\u\i}n and M.~A. Na{\u\i}mark.
\newblock The method of symmetric and {H}ermitian forms in the theory of the
  separation of the roots of algebraic equations.
\newblock {\em Linear and Multilinear Algebra}, 10(4):265--308, 1981.
\newblock Translated from the Russian by O. Boshko and J. L. Howland.

\bibitem{5}
A.~Kurosh.
\newblock {\em Higher algebra}.
\newblock ``Mir'', Moscow, 1988.
\newblock Translated from the Russian by George Yankovsky, Reprint of the 1972
  translation.

\bibitem{20}
S.~Lang.
\newblock {\em Algebra}, volume 211 of {\em Graduate Texts in Mathematics}.
\newblock Springer-Verlag, New York, third edition, 2002.

\bibitem{24}
J.~E. Marsden and M.~McCracken.
\newblock {\em The {H}opf bifurcation and its applications}.
\newblock Springer-Verlag, New York, 1976.
\newblock With contributions by P. Chernoff, G. Childs, S. Chow, J. R. Dorroh,
  J. Guckenheimer, L. Howard, N. Kopell, O. Lanford, J. Mallet-Paret, G. Oster,
  O. Ruiz, S. Schecter, D. Schmidt and S. Smale, Applied Mathematical Sciences,
  Vol. 19.

\bibitem{9}
G.~P{\'o}lya.
\newblock {\em Mathematics and plausible reasoning. {V}ol. {II}}.
\newblock Princeton University Press, Princeton, NJ, second edition, 1990.
\newblock Patterns of plausible inference.

\bibitem{8}
G.~Polya.
\newblock {\em How to solve it}.
\newblock Princeton Science Library. Princeton University Press, Princeton, NJ,
  2004.
\newblock A new aspect of mathematical method, Expanded version of the 1988
  edition, with a new foreword by John H. Conway.

\bibitem{10}
G.~P{\'o}lya.
\newblock {\em Mathematical discovery}.
\newblock John Wiley \& Sons Inc., New York, 1981.
\newblock On understanding, learning, and teaching problem solving, Reprint in
  one volume, With a foreword by Peter Hilton, With a bibliography by Gerald
  Alexanderson, With an index by Jean Pedersen.

\bibitem{11a}
G.~P{\'o}lya and G.~Szeg{\H{o}}.
\newblock {\em Problems and theorems in analysis. {I}}.
\newblock Classics in Mathematics. Springer-Verlag, Berlin, 1998.
\newblock Series, integral calculus, theory of functions, Translated from the
  German by Dorothee Aeppli, Reprint of the 1978 English translation.

\bibitem{11b}
G.~P{\'o}lya and G.~Szeg{\H{o}}.
\newblock {\em Problems and theorems in analysis. {II}}.
\newblock Classics in Mathematics. Springer-Verlag, Berlin, 1998.
\newblock Theory of functions, zeros, polynomials, determinants, number theory,
  geometry, Translated from the German by C. E. Billigheimer, Reprint of the
  1976 English translation.

\bibitem{2}
M.~M. Postnikov.
\newblock {\em Ustoichivye mnogochleny}.
\newblock ``Nauka'', Moscow, 1981.

\bibitem{16a}
T.~J. Stieltjes.
\newblock Recherches sur les fractions continues.
\newblock {\em Ann. Fac. Sci. Toulouse Math. (6)}, 4(1):Ji--Jiv, J1--J35, 1995.
\newblock Reprint of the 1894 original, With an introduction by Jean Cassinet.

\bibitem{16b}
T.~J. Stieltjes.
\newblock Recherches sur les fractions continues.
\newblock {\em Ann. Fac. Sci. Toulouse Math. (6)}, 4(2):Ji, J36--J75, 1995.
\newblock Reprint of the 1894 original.

\bibitem{16c}
T.~J. Stieltjes.
\newblock Recherches sur les fractions continues.
\newblock {\em Ann. Fac. Sci. Toulouse Math. (6)}, 4(3):J76--J122, 1995.
\newblock Reprint of Ann. Fac. Sci. Toulouse {\bf 8} (1894), J76--J122.

\bibitem{16d}
T.~J. Stieltjes.
\newblock Recherches sur les fractions continues.
\newblock {\em Ann. Fac. Sci. Toulouse Math. (6)}, 4(4):A5--A47, 1995.
\newblock Reprint of Ann. Fac. Sci. Toulouse {\bf 9} (1895), A5--A47.

\end{thebibliography}
\end{document}